\documentclass[12pt, a4paper]{amsart}

\usepackage{amsmath,amsfonts,amsthm,amssymb,indentfirst,epic}

\usepackage{hyperref}

\newtheorem{thm}{Theorem}[section]
\newtheorem{cor}[thm]{Corollary}
\newtheorem{lem}[thm]{Lemma}

\newtheorem{prop}[thm]{Proposition}

\newtheorem{defin}[thm]{Definition}

\newtheorem{remark}[thm]{Remark}
\newtheorem{remarks}[thm]{Remarks}

\def\R {{\mathbb R}}

\def\C {{\mathbb C}}
\def\N{{\mathbb N}}

\newcommand\ann{{\operatorname{ann}}}
\newcommand\core{{\operatorname{core}}}
\newcommand\GKdim{{\operatorname{GKdim}}}
\newcommand\Spec{{\operatorname{Spec}}}
\newcommand\Ht{{\operatorname{Ht}}}

\newcommand\kdim{{\operatorname{kdim}}}

\DeclareMathOperator{\End}{End} 
\DeclareMathOperator{\M}{M}

\newcommand{\ben}{\begin{enumerate}}
\newcommand{\een}{\end{enumerate}}
\newcommand{\bit}{\begin{itemize}}
\newcommand{\eit}{\end{itemize}}

\begin{document}
\title[]
{\large{Prime affine algebras of GK dimension two which are almost PI algebras}
 }
\author[]
{\small{Vered Moskowicz}}
\begin{abstract}
An almost PI algebra is a generalisation of a just infinite algebra which does not satisfy a polynomial identity.
An almost PI algebra has some nice properties: It is prime, has a countable cofinal subset of ideals and when satisfying ACC(semiprimes) (this happens, for example, when it is affine), it has only countably many height 1 primes.

Consider an affine prime Goldie non-simple non-PI $k$-algebra $R$ of GK dimension $<3$, where $k$ is an uncountable field. $R$ is an almost PI algebra.
We give some possible additional conditions which make such an algebra primitive, see Theorems \ref{algebraic images}, \ref{nilpotent images}, \ref{images} and \ref{images2}.
This gives a partial answer to Small's question: Let $R$ be an affine prime Noetherian $k$-algebra of GK dimension 2, where $k$ is any field. Does it follow that $R$ is PI or primitive?

We also show that the center of $R$ is a finite dimensional field extension of $k$ (Theorem \ref{fdfe}), and if, in addition, $k$ is algebraically closed, then $R$ is stably almost PI (Theorem \ref{stably almost}).
\end{abstract}

\maketitle
Remember that a just infinite algebra $R$ over a field $k$ is an algebra infinite dimensional over $k$ whose nonzero two sided ideals have finite codimension. A just infinite algebra can be PI or not.
A generalisation of a non-PI just infinite algebra was made by Rowen who suggested the notion of an almost PI algebra.
The notion of an almost PI algebra is also known to Small and probably to other researchers as well.

Let $R$ be an associative algebra with 1 over a field $k$.
$R$ is called \emph{an almost PI algebra} if $R$ is non-PI and $R/I$ is PI for each nonzero two sided ideal $I$ in $R$.
We will only deal with algebras over a field, although in some cases we can take a commutative Noetherian ring with 1 instead of a field.

In the first section we shall see some properties of almost PI algebras.

In the second section we shall try to answer the above mentioned question of Small:
Does it follow that an affine prime Noetherian $k$-algebra $R$ with GK dimension 2 is PI or primitive? 
Artin and Stafford in \cite{ast} have shown that this is true if $R$ is graded.
Bell has shown in \cite{bel} that this is false if $R$ is not Noetherian.
His example is of an affine prime algebra of GK dimension 2 over any field, which is non-PI and has a nonzero Jacobson radical (hence not primitive).
Smoktunowicz and Vishne's example in \cite{vish} also has a nonzero Jacobson radical.
Bell in \cite{bell} has shown that when $k$ is an uncountable field, an affine prime Goldie $k$-algebra with quadratic growth is PI or primitive. His result motivates the more general question: Does it follow that an affine prime Goldie $k$-algebra $R$ with GK dimension $<3$, where $k$ is an uncountable field, is PI or primitive?

We show that such an algebra $R$ is PI or primitive if, in addition, $k$ is algebraically closed, and there exists a non-algebraic element $x \in R$ such that $x+P \in R/P$ is algebraic over $k$, for every prime ideal $P$ with co-GK 1, see Theorem \ref{algebraic images}.
When $k$ is not algebraically closed, the existence of such a special non-algebraic element $x \in R$, may not guarantee that $R$ is primitive. Therefore, we add another assumption, see Theorem \ref{nilpotent images}, Theorem \ref{images} and Theorem \ref{images2}.
The proofs are based on ideas from the first section, from Farkas and Small's theorem (\cite[Theorem 2.2]{fs}) and from Bell's theorem (\cite[Theorem 1.1]{bell}). 

We also show that if $R$ is an affine prime Goldie non-simple non-PI $k$-algebra of GK dimension $<3$, where $k$ is an uncountable field, then its center is a finite dimensional field extension of $k$ (Theorem \ref{fdfe}) and if, in addition, $k$ is algebraically closed, then $R$ is stably almost PI (the $K$-algebra $R\otimes_k K$ is an almost PI algebra, for every field extension $K/k$. Theorem \ref{stably almost}).
The proofs are based on \cite{stably}.

In the third section we remark about just infinite algebras.

\section{Basic properties of almost PI algebras}

\subsection{Primality}
\begin{thm}\label{prime}
If $R$ is an almost PI algebra (not necessarily affine), then $R$ is prime.
\end{thm}
\begin{proof}
Otherwise, take $A,B$ nonzero ideals in $R$ such that $AB=0$. Let $f$ be a polynomial identity for $R/A$ and let $g$ be a polynomial identity for $R/B$. We can assume $f=f(X_1,\ldots,X_l)$ and $g=g(X_1,\ldots,X_l)$.
 Now, $\forall r_1,\ldots,r_l \in R$ we have $f(\bar{r_1},\ldots,\bar{r_l})=\bar{0}$, so $f(r_1,\ldots,r_l)+A=0+A$, hence $f(r_1,\ldots,r_l) \in A$.
Similarly, $\forall r_1,\ldots,r_l \in R$ $g(r_1,\ldots,r_l) \in B$.
 Therefore  $\forall r_1,\ldots,r_l \in R$
$(fg)(r_1,\ldots,r_l)=f(r_1,\ldots,r_l)g(r_1,\ldots,r_l) \in AB=0$, so $fg$ is a polynomial identity for $R$, which is impossible (by definition $R$ is non-PI).
\end{proof}

\begin{remark}
This proof also works for a non-PI just infinite algebra (but not for a PI just infinite algebra).
\end{remark}

The next is a very simple observation.
\begin{cor}\label{1.cor.simple}
Let $R$ be an almost PI algebra. \bit
\item [(i)] If $R$ is Artinian, then it is simple.
\item [(ii)] If $R$ is Goldie and algebraic, then it is simple.
\eit
\end{cor}

\begin{proof}
\bit
\item [(i)] Any prime Artinian ring is simple \cite[2.3.9]{ro1}.
\item [(ii)] Any algebra which is prime Goldie and algebraic is simple. Indeed, in a prime Goldie algebra each nonzero two sided ideal contains a regular element \cite[Proposition 3.2.13]{ro1}. But this regular element is algebraic, hence invertible.
\eit
\end{proof}

There is a criterion for an affine infinite dimensional algebra to be just infinite: every nonzero prime ideal has finite codimension, see \cite[Corollary 1]{prime}.
Instead of affinity we will assume ACC(ideals) and get the following criterion, which can be thought of as a corollary to Theorem \ref{prime}.

\begin{cor}\label{criterion}
Let $R$ be a weakly Noetherian non-PI algebra. Assume that for every nonzero prime ideal $P$, $R/P$ is PI.
Then $R$ is an almost PI algebra.
\end{cor}

\begin{proof}
Otherwise, there exists a nonzero ideal $\tilde{I}$ with $R/\tilde{I}$ non-PI.
So $T:=\{0 \neq I\triangleleft R | R/I \ \text{is non-PI} \} \neq \emptyset$.
Let $I_1 \subseteq I_2 \subseteq I_3 \subseteq \ldots$ be an ascending chain of ideals in $T$. Since $R$ satisfies ACC(ideals), it follows that there exists $m \geq 1$ such that $I_m=I_{m+1}=I_{m+2}=\ldots$. So $I_m \in T$ is an upper bound for our chain. Hence Zorn's Lemma implies that $T$ contains a maximal element, call it $I$.

$R/I$ is an almost PI algebra; Indeed, $I \in T$ implies that $R/I$ is non-PI, and maximality of $I$ implies that for every $I \subsetneq J$, $R/J$ is PI, so $(R/I)/(J/I)\cong R/J$ is PI.

Therefore, from Theorem \ref{prime}, $R/I$ is prime. Hence $I$ is a (nonzero) prime ideal of $R$, so by assumption $R/I$ is PI.
We reached a contradiction, so there exists no such nonzero ideal $\tilde{I}$ with $R/\tilde{I}$ non-PI.
\end{proof}

Of course, a similar criterion can be given for a just infinite algebra. Observe that a just infinite algebra $A$ is actually weakly Noetherian.
Indeed, let $I_1 \subsetneq I_2 \subsetneq I_3 \subsetneq \ldots$ be an ascending chain of proper ideals in $A$. 
Let $I:= \cup I_i$. 
$A/I_1\stackrel{onto}{\rightarrow} A/I_2 \stackrel{onto}{\rightarrow} \ldots \stackrel{onto}{\rightarrow} A/I$, implies $dim_k(A/I)<\ldots <dim_k(A/I_2)<dim_k(A/I_1)<\infty$, which is impossible.

\begin{prop}\label{criterion2}
Let $R$ be a weakly Noetherian infinite dimensional $k$-algebra. Assume that for every nonzero prime ideal $P$, $R/P$ is finite dimensional over $k$.
Then $R$ is a just infinite algebra.
\end{prop}

\begin{proof}
Just replace in the proof of Corollary \ref{criterion}: non-PI by infinite dimensional, PI by finite dimensional, almost PI by just infinite, and Theorem \ref{prime} by \cite[Theorem 1]{prime} (which says that a just infinite algebra is prime, whether affine or not).
\end{proof}

Next, we mimick an observation of Small (\cite[Proposition 3.2]{PassTem} or \cite[Lemma 2.1]{fs}), but again instead of affinity we assume ACC(ideals).

\begin{lem}
Let $R$ be a weakly Noetherian non-PI $k$-algebra. Then $R$ has a homomorphic image $\bar{R}$ which is prime and almost PI. 
\end{lem}
\begin{proof}
Define $T:=\{I\triangleleft R | R/I \ \text{is non-PI} \}$. Since $R\cong R/0$ is non-PI, $T \neq \emptyset$.
Let $I_1 \subseteq I_2 \subseteq I_3 \subseteq \ldots$ be an ascending chain of ideals in $T$. Since $R$ satisfies ACC(ideals), it follows that there exists $m \geq 1$ such that $I_m=I_{m+1}=I_{m+2}=\ldots$. So $I_m \in T$ is an upper bound for our chain. Hence Zorn's Lemma implies that $T$ contains a maximal element, call it $I$.

$R/I$ is an almost PI algebra; Indeed, $I \in T$ implies that $R/I$ is non-PI, and maximality of $I$ implies that for every $I \subsetneq J$, $R/J$ is PI, so $(R/I)/(J/I)\cong R/J$ is PI.

Now there are two ways to see that $\bar{R}=R/I$ is prime:
First, just use Theorem \ref{prime}.
Second, use the following general claim:
Let $S$ be a non-PI algebra. If $A$ and $B$ are (nonzero) ideals of $S$ such that $S/A$ and $S/B$ are PI, then $AB \neq 0$: Otherwise, $AB=0$. Let $f$ be a polynomial identity for $S/A$ and let $g$ be a polynomial identity for $S/B$. We can assume $f=f(X_1,\ldots,X_l)$ and $g=g(X_1,\ldots,X_l)$.
Now, $\forall s_1,\ldots,s_l \in S$ we have $f(\bar{s_1},\ldots,\bar{s_l})=\bar{0}$, so $f(s_1,\ldots,s_l)+A=0+A$, hence $f(s_1,\ldots,s_l) \in A$.
Similarly, $\forall s_1,\ldots,s_l \in S$ $g(s_1,\ldots,s_l) \in B$.
Therefore, $\forall r_1,\ldots,s_l \in S$, $(fg)(s_1,\ldots,s_l)=f(s_1,\ldots,s_l)g(s_1,\ldots,s_l) \in AB=0$, so $fg$ is a polynomial identity for $S$, which is impossible (by assumption $S$ is non-PI).

In our case, $\bar{R}$ is a non-PI algebra (it is almost PI). Take any nonzero ideals $\bar{A}$ and $\bar{B}$ in $\bar{R}$. $\bar{R}/\bar{A}$ and $\bar{R}/\bar{B}$ are PI, so from the above claim $\bar{A}\bar{B} \neq \bar{0}$. So $\bar{R}$ is prime.

\end{proof}

Small's observation shows that an affine just infinite algebra is prime (however, any just infinite algebra is prime \cite[Theorem 1]{prime}) and our observation shows that a weakly Noetherian almost PI algebra is prime (however, any almost PI algebra is prime as was seen above in Theorem \ref{prime}).

\subsection{Primitive ideals}

In a PI ring every primitive ideal is maximal. Almost the same is true in our case:
\begin{prop}\label{maximal}
If $R$ is an almost PI algebra, then every nonzero primitive ideal is maximal.
\end{prop}

\begin{proof}
Take a nonzero primitive ideal $P$ in $R$.
$R/P$ is primitive PI so by Kaplansky's theorem it is simple.
So $P$ is maximal.
\end{proof}

The next Proposition says that in an \textit{affine} almost PI algebra, every nonzero primitive ideal has finite codimension. This property will be used in Lemma 3.5.

\begin{prop}\label{1.7prop}
Let $R$ be an almost PI algebra and let $P$ be a nonzero prime ideal in $R$. Then:\bit
\item [(i)] $\dim_k(R/P)<\infty \Rightarrow P \text{\ is primitive (=maximal)}$.
\item [(ii)] If $R$ is, in addition, affine: $\dim_k(R/P)<\infty \Leftrightarrow P \text{\ is primitive}$.
\eit
\end{prop}

\begin{proof}
\bit
\item [(i)] $R/P$ is prime and $\dim_k(R/P)<\infty$, hence $R/P$ is simple (\cite[Corollary 1.6.30]{ro80}), so $P$ is maximal.
\item [(ii)] ($\Leftarrow$): $R/P$ is primitive PI so by Kaplansky's theorem it is Artinian. So $R/P$ is PI affine and Artinian, hence finite dimensional \cite[13.10.3(iv)]{mr}.
\eit
\end{proof}

\begin{lem}\label{1.8lem}
Let $R$ be an affine almost PI algebra. If:
\bit
\item[(i)] For each nonzero prime ideal $P$, $\dim_k(R/P)<\infty$.

or more generally: \item[(ii)] For each nonzero prime ideal $P$, $R/P$ is algebraic over $k$.
\eit
Then $R$ is non-PI just infinite.
\end{lem}
\begin{proof}
\bit
\item [(i)] This follows from \cite[Corollary 1]{prime} which says that if in an affine infinite dimensional algebra each nonzero prime ideal has finite codimension, then the algebra is just infinite (of course our algebra is infinite dimensional, since it is non-PI).
\item [(ii)] For each nonzero prime ideal $P$, $R/P$ is an affine algebraic PI algebra, hence finite dimensional (Kurosch problem for PI, see \cite[13.8.9(2)]{mr}).
\eit
\end{proof}

Form Proposition \ref{1.7prop} and Lemma \ref{1.8lem} we get the following corollaries.

\begin{cor}\label{1.9cor}
Let $R$ be an affine almost PI algebra.
If every nonzero prime ideal is primitive, then $R$ is non-PI just infinite.
\end{cor}

\begin{proof}
From Proposition \ref{1.7prop}(ii) every nonzero prime ideal has finite codimension. Now use Lemma \ref{1.8lem}(i).
\end{proof}

\begin{cor}\label{1.classical.krull}
Let $R$ be an affine almost PI algebra and assume that it is not simple.

$R$ is non-PI just infinite  $\Leftrightarrow$ $R$ has classical Krull dimension 1.
\end{cor}

\begin{proof}
($\Rightarrow$) Any just infinite algebra (whether PI or not, whether affine or not) has classical Krull dimension 1; Indeed, each nonzero prime ideal has finite codimension, hence is maximal (see Proposition \ref{1.7prop}(i)). Also any just infinite algebra is prime (\cite[Theorem 1]{prime}). Actually, since now we are only interested in a non-PI just infinite (affine) algebra, we have seen above in Remark 1.2 that such an algebra is prime. Hence it is clear that $R$ has classical Krull dimension 1.

($\Leftarrow$) In view of Lemma \ref{1.8lem}(i) it is enough to show that each nonzero prime ideal has finite codimension. Let $P$ be a nonzero prime ideal. Since $R$ has classical Krull dimension 1, $P$ must be maximal (otherwise, there exists a maximal ideal $Q$ strictly containing $P$: $0\subsetneq P \subsetneq Q$, a contradiction). Now, Proposition \ref{1.7prop}(ii) implies $\dim_k(R/P)<\infty$.
(Another proof for ($\Leftarrow$): Since $R$ has classical Krull dimension 1, every nonzero prime ideal is maximal. Hence, from corollary \ref{1.9cor} $R$ is non-PI just infinite).
\end{proof}

The following proposition shows that over an uncountable field $k$, an affine (non-simple) almost PI algebra has classical Krull dimension 1, if it has only countably many primitive ideals and if it has finite classical Krull dimension.

\begin{remark}\label{rowen}
In Rowen's paper \cite{rosem}:
The conclusion of Proposition 3.2 that $R$ is primitive is true even if $R$ only has a countable separating set (which satisfies the additional restriction (4) of that paper) for its nonzero primitive ideals instead of for its nonzero ideals. This can be seen directly from the proof there.
\end{remark}

\begin{prop}\label{classical}
Let $R$ be an affine non-simple almost PI algebra over an uncountable field having only countably many primitive ideals. If $R$ has finite classical Krull dimension, then $R$ has classical Krull dimension 1.
In other words, there is no such algebra $R$ with $2 \leq \kdim(R)< \infty$, where $\kdim(R)$ is the classical Krull dimension of $R$. Moreover, $R$ is non-PI just infinite.
\end{prop}

\begin{proof}
By induction on $\kdim(R)< \infty$.

$\kdim(R)=2$: We shall see that this is impossible.
Otherwise, let $R$ be such an algebra. Denote all height 1 primes in $R$ by $\{P_\alpha\}_{\alpha\in A}$ and all height 2 primes in $R$ by $\{Q_\beta\}_{\beta\in B}$. So $\Spec(R)=\{0,P_\alpha,Q_\beta\}$.
$|B|\leq\aleph_0$ from our assumption that there are only countably many primitive ideals. Indeed, each height 2 prime must be maximal, since the classical Krull dimension of $R$ is 2.
(Notice: $|A|\leq\aleph_0$ from Proposition \ref{aa}).
We claim that each height 1 prime must be maximal: Otherwise, there exists a height 1 prime $P_\alpha$ which is not maximal. So $P_\alpha$ is strictly contained in some height 2 prime. Let $\{Q_i\}_{i\in I}$ be all those height 2 primes which contain $P_\alpha$ ($I\subseteq B$).
We shall now see that $P_\alpha$ must be primitive, hence by Proposition \ref{maximal} maximal, a contradiction to our choice of $P_\alpha$ not being maximal.
\bit
\item If $|I|<\infty$: $\bar{0}=J(R/P_\alpha)=\cap_{i\in I} (Q_i/P_\alpha)$ where $R/P_\alpha$ is semiprimitive by primality and the Razmyslov-Kemer-Braun theorem, and since $R/P_\alpha$ is prime (of course in a prime ring finite intersections of nonzero ideals are nonzero), it follows that $\bar{0}$ must be a primitive ideal of $R/P_\alpha$. \item If $|I|=\aleph_0$: actually in this case we shall see a claim also applicable to $|I|<\infty$, so it is possible not to separate to two cases.
$\{Q_i/P_\alpha\}_{i\in I}$ are all the nonzero prime ideals of $R/P_\alpha$, and each is maximal in $R/P_\alpha$ (since each $Q_i$ is maximal in $R$). So $\{Q_i/P_\alpha\}_{i\in I}$ are all the nonzero primitive ideals of $R/P_\alpha$. Therefore, taking $\bar{0}\neq\bar{y_i}\in (Q_i/P_\alpha)$ yields a countable separating set for the nonzero primitive ideals of $R/P_\alpha$. Since $R/P_\alpha$ is PI, by Rowen's theorem $Z(R/P_\alpha)\cap\bar{R}\bar{y_i}\bar{R}\neq \bar{0}$, hence we may take $\bar{0}\neq \bar{x_i}\in Z(R/P_\alpha)\cap\bar{R}\bar{y_i}\bar{R}$ and get a countable separating set for the nonzero primitive ideals of $R/P_\alpha$ with elements regular (any nonzero element of the center of a prime ring is regular) and commuting with each other. Thus from Theorem 3.3 of \cite{rosem} this separating set satisfies (4) in that paper. Hence we can use Proposition 3.2 there (with a little generalisation, as mentioned in Remark \ref{rowen}) and get that $R/P_\alpha$ is primitive. \eit
Therefore, each height 1 prime must be maximal. But then obviously primes of height 2 cannot exist, hence $\kdim(R)=1$ (remember that we have assumed that $R$ is not simple, so $\kdim(R) \geq 1$).

$\kdim(R)=n$: We shall see that this is impossible.
Otherwise, let $R$ be such an algebra. Denote all height 1 primes in $R$ by $\{P_\alpha\}_{\alpha\in A}$, all height 2 primes by $\{Q_\beta\}_{\beta\in B}$, \ldots , all height $n-1$ primes by $\{F_\nu\}_{\nu\in M}$, all height $n$ primes by $\{G_\nu\}_{\nu\in N}$. So $\Spec(R)=\{0,P_\alpha,Q_\beta,\ldots,F_\mu,G_\nu\}$.
$|N|\leq\aleph_0$ from our assumption that there are only countably many primitive ideals. Indeed, each height n prime must be maximal, since the classical Krull dimension of $R$ is $n$.
(Notice: $|A|\leq\aleph_0$ from Proposition \ref{aa}). Be careful: we only know that $|A|=|N|\leq\aleph_0$, but the cardinality of the height 2 primes until the height $n-1$ primes is not known.
We claim that each height $n-1$ prime must be maximal: Otherwise, there exists a height $n-1$ prime $F_\mu$ which is not maximal. So $F_\mu$ is strictly contained in some height $n$ prime. Let $\{G_i\}_{i\in I}$ be all those height $n$ primes which contain $F_\mu$ ($I\subseteq N$).
We shall now see that $F_\mu$ must be primitive, hence by Proposition \ref{maximal} maximal, a contradiction to our choice of $F_\mu$ not being maximal.
\bit
\item If $|I|<\infty$: $\bar{0}=J(R/F_\mu)=\cap_{i\in I} (G_i/F_\mu)$ where $R/F_\mu$ is semiprimitive by primality and the Razmyslov-Kemer-Braun theorem, and since $R/F_\mu$ is prime (of course in a prime ring finite intersections of nonzero ideals are nonzero), it follows that $\bar{0}$ must be a primitive ideal of $R/F_\mu$.
\item If $|I|=\aleph_0$: actually in this case we shall see a claim also applicable to $|I|<\infty$, so it is possible not to separate to two cases.
$\{G_i/F_\mu\}_{i\in I}$ are all the nonzero prime ideals of $R/F_\mu$, each is maximal in $R/F_\mu$ (since each $G_i$ is maximal in $R$). So $\{G_i/F_\mu\}_{i\in I}$ are all the nonzero primitive ideals of $R/F_\mu$. Therefore, taking $\bar{0}\neq\bar{y_i}\in (G_i/F_\mu)$ yields a countable separating set for the nonzero primitive ideals of $R/F_\mu$. Since $R/F_\mu$ is PI, by Rowen's theorem $Z(R/F_\mu)\cap\bar{R}\bar{y_i}\bar{R}\neq \bar{0}$, hence we may take $\bar{0}\neq \bar{x_i}\in Z(R/F_\mu)\cap\bar{R}\bar{y_i}\bar{R}$ and get a countable separating set for the nonzero primitive ideals of $R/F_\mu$ with elements regular (any nonzero element of the center of a prime ring is regular) and commuting with each other. Thus from Theorem 3.3 of \cite{rosem} this separating set satisfies (4) in that paper. Hence we can use Proposition 3.2 there (with a little generalisation, as mentioned in Remark \ref{rowen}) and get that $R/F_\mu$ is primitive. \eit
Therefore each height $n-1$ prime must be maximal. But then obviously primes of height $n$ cannot exist, hence  $\kdim(R)\leq n-1$, so by the induction hypothesis $\kdim(R)=1$.

Finally, $R$ is non-PI just infinite from Corollary \ref{1.classical.krull}.
\end{proof}

\begin{remark}
If one reads carefully the proof of Proposition \ref{classical}, one sees that instead of "affine over an uncountable field", we could assume "the usual dimension-cardinality hypothesis".
This is because:
\bit
\item "(Notice: $|A|\leq\aleph_0$ from Proposition \ref{aa})": In Proposition \ref{aa} affinity is needed, but the countability of the height 1 primes is not used in the proof.
\item "We shall now see that $P_\alpha$ must be primitive, hence by Proposition \ref{maximal} maximal": In Proposition \ref{maximal} affinity is not needed.
\item In Proposition 3.2 and in Theorem 3.3 of \cite{rosem} affiniy is not needed, only $dim_k(R)<|k|$.
\eit
\end{remark}

\subsection{Nullstellensatz}
Another property of an affine almost PI algebra can be obtained by using the Razmyslov-Kemer-Braun theorem \cite[Theorem 2.57]{belov}, which says that every affine PI algebra has a nilpotent Jacobson radical.

A Jacobson ring is a ring in which every prime factor ring is semiprimitive.
A ring $R$ has the \emph{radical property} if the Jacobson radical of each factor ring of $R$ is nil.

Each Jacobson ring has the radical property, see \cite[9.1.2]{mr}.

\begin{prop}\label{1.15}
Let $R$ be an affine almost PI algebra. For every nonzero prime ideal $P$, $R/P$ is semiprimitive.
Therefore, if $R$ is semiprimitive, then it is a \emph{Jacobson ring}.
\end{prop}

\begin{proof}
Let $P$ be a nonzero prime ideal of $R$. $R/P$ is an affine PI algebra, so by Razmyslov-Kemer-Braun theorem $J(R/P)$ is nilpotent. Primality of $R/P$ implies $J(R/P)=0$.
The second part is obvious (remember that we have seen in Theorem \ref{prime} that an almost PI algebra is prime).
\end{proof}

An algebra $R$ over a field $k$ has the \emph{endomorphism property over} $k$ if, for each simple $R$-module $M$, $\End_R(M)$ is algebraic over $k$.

An algebra $R$ over a field $k$ \emph{satisfies the Nullstellensatz over} $k$ if $R$ has both the radical property and the endomorphism property.

From Proposition \ref{1.7prop}(ii) we get the following proposition.
\begin{prop}\label{1.16}
If $R$ is an affine almost PI algebra which is not primitive, then each simple $R$-module is finite dimensional over $k$. Therefore $R$ has the endomorphism property.
\end{prop}

\begin{proof}
Take a maximal left ideal $L$ in $R$. Since we have assumed that $R$ is not primitive, $\core(L)\neq0$. Denote $P=\core(L)$. $P$ is a nonzero primitive ideal in $R$ therefore, from Proposition \ref{1.7prop}(ii), $P$ has finite codimension. From $P=\core(L)\subseteq L$ we see $L/P\subseteq R/P$, hence $(R/P)/(L/P)\cong R/L$ is finite dimensional over $k$.
Since any simple $R$-module is isomorphic to $R/L$ where $L$ is a maximal left ideal in $R$, we are done.
Finally, $R$ has the endomorphism property: Let $M$ be a simple $R$-module. $\End_R(M)\subseteq \End_k(M)$ and $\End_k(M)$ is finite dimensional over $k$ (since $M$ is finite dimensional over $k$), so obviously $\End_k(M)$ is algebraic over $k$. Hence $\End_R(M)$ is also algebraic over $k$.
\end{proof}

Proposition \ref{1.16} tells us that if an affine almost PI algebra $R$ has a maximal left ideal $L$ with $\dim_k(R/L)=\infty$, then $R$ is primitive.

From Proposition \ref{1.15} and Proposition \ref{1.16} we get the following corollary.
When $k$ is uncountable, it is already known that any countably generated $k$-algebra satisfies the Nullstellensatz, see \cite[9.1.8]{mr}.
Hence, the next corollary adds something new only when $|k| \leq \aleph_0$.

\begin{cor}
If $R$ is an affine almost PI algebra which is semiprimitive but not primitive, then $R$ satisfies the Nullstellensatz.
More generally, it is enough to demand that $J(R)$ is nil (instead of $J(R)=0$).
\end{cor}

\begin{proof}
Obvious from Proposition \ref{1.15} and Proposition \ref{1.16}.
Also the more general claim is obvious.
\end{proof}

\subsection{Countable cofinal subset of ideals}
\begin{defin}
Fisher and Snider \cite{von-Neumann} have defined for a prime ring a \emph{countable cofinal subset of ideals} as a countable collection of nonzero ideals, such that any nonzero ideal contains an ideal in that collection.
\end{defin}

\begin{defin}
By a \emph{separating set $S$ for the nonzero ideals of $R$} we mean a set $S$ $(0\notin S)$ such that for each nonzero ideal $I$ in $R$, $I\cap S\neq\emptyset$.
\end{defin}

The following property of an almost PI algebra (not necessarily affine) will be used many times in this paper.
\begin{prop}\label{ccsubset}
If $R$ is an almost PI algebra, then $R$ has a countable cofinal subset of ideals. Therefore $R$ has a countable separating set for its nonzero ideals.
\end{prop}

\begin{proof}
Define for every $n,d \in \N$  $S_{n,d}:=\left\langle ((s_n)^d)(R) \right\rangle$ the two sided ideal generated by the set $((s_n)^d)(R)=\{((s_n)^d)(r_1,\ldots,r_n)|r_i \in R\}$ where $s_n$ is the standard polynomial of degree $n$. Each $S_{n,d}$ is nonzero, since $R$ is non-PI.

Now take any $0\neq I\triangleleft R$. $R/I$ is PI, so it satisfies some power of a standard polynomial \cite[Theorem 13.4.8]{mr}, assume it to be $(s_n)^d$ for some $n,d \in \N$.

$\forall \bar{r_i} \in R/I \ ((s_n)^d)(\bar{r_1},\ldots,\bar{r_n})=\bar{0}$.

$\forall r_i \in R \ ((s_n)^d)(r_1,\ldots,r_n)+I=0+I$,

$\forall r_i \in R \ ((s_n)^d)(r_1,\ldots,r_n)\in I$,
hence $((s_n)^d)(R)\subseteq I$, implying $S_{n,d} \subseteq I$.

That $R$ has a countable separating set for its nonzero ideals follows at once: just choose one nonzero element $y_{n,d}$ from each $S_{n,d}$, therefore $\{y_{n,d}\}_{n,d \in \N}$ is a countable separating set.
\end{proof}

\begin{remarks}
\bit
\item [(i)] In an \textit{affine} almost PI algebra there exists another countable cofinal subset of ideals, namely $\{T_n(R)\}_{n \in  \N}$, where $T_n(R)$ is the two sided ideal generated by the set $t_n(R)=\{f(r_1,\ldots,r_m)|f \in id(\M_n(k)),r_i \in R\}$.
Indeed; First, each $T_n(R)$ is nonzero, since $R$ is non-PI.
Second, take any $0\neq I\triangleleft R$. $R/I$ is affine PI, so it satisfies all polynomial identities of $\M_n(k)$ for some $n$ (since $J(R/I)$ is nilpotent and \cite[page 216]{ro80} or \cite[Theorem 6.3.16]{ro2}):

$\forall f=f(X_1,\ldots,X_m)\in id(\M_n(k)) \ \ \forall \bar{r_i} \in R/I \ f(\bar{r_1},\ldots,\bar{r_m})=\bar{0}$.
Then, for every $f=f(X_1,\ldots,X_m)\in id(\M_n(k))$ and every $r_i \in R$, we have that $f(r_1,\ldots,r_m)+I=0+I$, so, for every $f=f(X_1,\ldots,X_m)\in id(\M_n(k))$ and every $r_i \in R$, we must have $f(r_1),\ldots,r_m)\in I$,
which means that $t_n(R)\subseteq I$, hence $T_n(R)\subseteq I$.

\item [(ii)] If a prime Goldie ring $R$ has a separaring set $S$, then we can assume each element in $S$ is regular. This is because we can replace $s\in S$ by a regular element $x\in RsR$ ($RsR$ contains a regular element since it is a nonzero ideal in a prime Goldie ring). Hence any Goldie almost PI algebra has a countable separating set of regular elements.
\eit
\end{remarks}

Remember that if $R$ is any ring which satisfies ACC(semiprimes), then: \bit
\item [(a)] every semiprime ideal is a finite intersection of prime ideals,
\item [(b)] there are only finitely many prime ideals minimal over any ideal. \eit
(see \cite[Proposition 5.2.3 and Corollary 5.2.4]{pro}). 

Using (b) and an idea from \cite[Proposition 2.4]{rosem} gives the following proposition.
It is very important for the proofs of Theorem \ref{finitecoGK1}, Theorem \ref{algebraic images} and Theorem \ref{nilpotent images}.

\begin{prop}\label{apACC}
If $R$ is an almost PI algebra satisfying ACC(semiprimes), then $R$ has only countably many height 1 primes.
\end{prop}

\begin{proof}
From proposition \ref{ccsubset}, $R$ has a countable separating set for its nonzero ideals, denote it $S=\{s_i\}_{i \in \N}$.
For each $Rs_iR$, let $P_{i1},\ldots,P_{it(i)}$ be those prime ideals of $R$ minimal over $Rs_iR$. Indeed, there are only finitely many such primes since we have assumed that $R$ satisfies ACC(semiprimes). Let $\Phi=\{P_{iu}| i \in \N, 1\leq u \leq t(i), \Ht(P_{iu})=1\}$. $\Phi$ is of course a countable set of height 1 primes. We will see now that every height 1 prime belongs to $\Phi$: Take any height 1 prime $Q$. Since $S$ is a separating set, $s_j \in Q$ for some $j \in \N$. Therefore $Rs_jR\subseteq Q$. $Q$ must be minimal over $Rs_jR$ from considerations of height ($Q$ is of height 1). Thus $Q \in \Phi$.
\end{proof}

\begin{remarks}\label{1.remarks}
(i) In the proof of Proposition \ref{apACC} it is enough to have a countable separating set for the nonzero height 1 prime ideals of $R$; However, from Proposition \ref{ccsubset}, $R$ has a countable separating set for all its nonzero ideals.

(ii) If $R$ is an almost PI algebra satisfying ACC(semiprimes), then any infinite (necessarily countable, from Proposition \ref{apACC}) intersection of height 1 primes is zero; Otherwise, there are height 1 primes $\{Q_n\}_{n \in \N}$ with $\cap Q_n \neq0$. Since $\cap Q_n$ is a semiprime ideal, it must be a finite intersection of prime ideals, as was mentioned in (a) above. Thus there exist $P_1,\ldots,P_m$ nonzero primes such that $P_1\cap \cdots \cap P_m= \cap Q_n$. But this is impossible because then $\forall j \in \N$: $P_1 \cdots P_m\subseteq P_1\cap \cdots \cap P_m= \cap Q_n\subseteq Q_j$, which implies $\{P_i\}_{1\leq i \leq m} \supseteq\{Q_n\}$, a contradiction.

Actually, one sees that in any prime algebra (not necessarily almost PI) which satisfies ACC(semiprimes), any infinite (not necessarily countable) intersection of height 1 primes is zero.
\end{remarks}

From the above remarks we can generalize Proposition \ref{apACC}. Instead of an almost PI algebra we can take a prime non-PI algebra such that for every nonzero prime ideal, its factor ring is PI (or even more generally: for every nonzero \emph{height} 1 prime ideal, its factor ring is PI. This is a generalisation if that algebra does not satisfy DCC(primes)).

\begin{prop}\label{aapACC}
Let $R$ be a prime non-PI ring such that:
\ben
\item[(i)] For every nonzero prime ideal $P$, $R/P$ is PI.
\item[(ii)] $R$ satisfies ACC(semiprimes).
\een
Then $R$ has only countably many height 1 primes.

More generally, 
Instead of (i) we can assume: (i') For every nonzero height 1 prime $P$, $R/P$ is PI,
and still $R$ has only countably many height 1 primes.
\end{prop}

\begin{proof}
(i) implies that $R$ has a countable separating set for its nonzero prime ideals. Indeed, take a nonzero prime ideal $P$. $R/P$ is PI, so it satisfies some power of a standard polynomial, assume it to be $(s_n)^d$ for some $n,d \in \N$. Hence $S_{n,d} \subseteq P$, where $S_{n,d}:=\left\langle (s_n)^d(R) \right\rangle$ the two sided ideal generated by the set $(s_n)^d(R)=\{(s_n)^d(r_1,\ldots,r_n)|r_i \in R\}$.
Now just choose one nonzero element $y_{n,d}$ from each $S_{n,d}$ ($R$ is non-PI, so $S_{n,d}\neq 0$).
Therefore $\{y_{n,d}\}_{n,d \in \N}$ is a countable separating set for the nonzero prime ideals of $R$, rename it $S=\{s_i\}_{i \in \N}$.
Reading carefully the proof of Proposition \ref{apACC} reveals that the same proof works here also, only now $S=\{s_i\}_{i \in \N}$ is a countable separating set for the nonzero prime ideals.

When assuming (i') instead of (i): Just add "height 1" where requiered. 
\end{proof}

Now we will consider affine almost PI algebras.

\begin{lem}\label{1.ACCsemiprimes}
If $R$ is an affine almost PI algebra, then $R$ satisfies ACC(semiprimes).
\end{lem}
\begin{proof}
Let $0\neq \tilde{P_1}\subset \tilde{P_2}\subset \tilde{P_3}\subset \cdots$ be an ascending chain of semiprime ideals in $R$.
Then in $R/\tilde{P_1}$ we have $\tilde{P_2}/\tilde{P_1}\subset \tilde{P_3}/\tilde{P_1}\subset \cdots$, which is an ascending chain of semiprime ideals.
This is impossible, since $R/\tilde{P_1}$ is an affine PI algebra, and as such it must satisfy the ascending chain condition on semiprime ideals \cite[Corollary 5.2.2]{pro}.
\end{proof}

The following Proposition is what will actually be used in the proofs of Theorem \ref{finitecoGK1}, Theorem \ref{algebraic images}, Theorem \ref{nilpotent images} and Theorem \ref{images}.

\begin{prop}\label{aa}
If $R$ is an affine almost PI algebra, then $R$ has only countably many height 1 primes.
\end{prop}
\begin{proof}
From Lemma \ref{1.ACCsemiprimes} $R$ satisfies ACC(semiprimes). Now use Proposition \ref{apACC}.
\end{proof}

\subsection{von Neumann regularity}
Notice that affinity is not needed in the following theorem.
\begin{thm}\label{vnr}
If an almost PI algebra $R$ is von Neumann regular, then it is primitive.
\end{thm}

\begin{proof}
$R$ is an almost PI algebra, so it is prime by Theorem 1.1 and has a countable cofinal subset of ideals by Proposition \ref{ccsubset}. Now use Corollary 1.2 in \cite{von-Neumann}.
\end{proof}

\subsection{Algebraicity}
Remember that von Neumann regularity implies semiprimitivity \cite[2.11.19(ii)]{ro1}.
It is then interesting to see what happens if instead of von Neumann regularity we only assume semiprimitivity.
Since now it is a more general question, it seems reasonable to add another assumption to compensate for what was lost.
Again \cite{von-Neumann} will help us.

The idea of the proof of the following theorem is taken from the proof of Farkas and Small's Theorem \cite[Theorem 2.2]{fs}. 

\begin{thm}\label{semialg}
Let $k$ be any field and let $R$ be an algebraic almost PI $k$-algebra. Then $J(R)\neq 0$ or $R$ is primitive.
\end{thm}

\begin{proof}
Assume that $J(R)=0$.
Therefore, each nil left ideal is zero. Hence every nonzero left ideal is not nil, so it has a non nilpotent element. By a remark of Rowen \cite[page 324]{rosem}, every nonzero left ideal of $R$ contains an idempotent different from 0,1. $R$ is prime and has a countable cofinal subset of ideals.
Then Theorem 1.1 in \cite{von-Neumann} implies that $R$ is primitive.
\end{proof}

\begin{cor}\label{1.29}
Let $k$ be any field and let $R$ be a semiprimitive $k$-algebra.
Assume that the following two conditions are satisfied:
\bit
\item [(i)] $R/I$ is PI, for every nonzero ideal $I$ in $R$.
\item [(ii)] $R$ is algebraic over $k$,
\eit
Then $R$ is PI or primitive. 
\end{cor}

Notice that one may consider Corollary \ref{1.29} interesting only when $R$ is not prime Goldie (or when $R$ is not a domain). Indeed, an algebraic prime Goldie algebra is simple- an explanation is given in the beginning of section 2 (an algebraic domain is a division algebra).

\begin{proof}
Assume that $R$ is non-PI.
Then, from (i), $R$ is an almost PI algebra. So Theorem \ref{semialg} implies that $R$ is primitive.
\end{proof}

We have already seen the following theorem before (Lemma \ref{1.8lem}). It is an immediate consequence of a solution, for PI rings, to Kurosch problem:
\begin{thm}\label{1.27thm}
If an affine almost PI $k$-algebra $R$ is algebraic, then it is non-PI just infinite.  
\end{thm}

\begin{proof}
Let $I$ be a nonzero ideal of $R$. $R/I$ is an algebraic affine PI algebra, hence finite dimensional \cite[13.8.9(ii)]{mr}.
\end{proof}

In Theorem \ref{1.27thm} we have taken an algebraic affine almost PI algebra. Now we want to see what may make an affine almost PI algebra algebraic. 

\begin{prop}
Let $k$ be an uncountable field and let $R$ be an affine almost PI $k$-algebra.
Assume that $R$ has only countably many primitive ideals.
Then $R$ is algebraic over $k$ or primitive.
\end{prop}

\begin{proof}
The proof is actually Farkas and Small's proof \cite[Theorem 2.2]{fs}.
We only give a sketch of it, since in section 2 and in section 3, detailed similar proofs will be given.
Assume that $R$ is not primitive.
We must show that $R$ is algebraic over $k$.
Otherwise, there exists a non-algebraic element $x \in R$.
Define two subsets of $k$, $A$ and $B$, as in Theorem \ref{fs}.
Notice that:\bit
\item In Theorem \ref{fs}, non-PI implies that there are only countably many primitive ideals, while here we assume that there are only countably many primitive ideals.
\item There as well as here, the countability of $B$ follows from the finite codimensionality of the nonzero primitive ideals. There the nonzero primitive ideals has finite codimension since the algebra is just infinite, while here the nonzero primitive ideals has finite codimension since the algebra is affine almost PI, see Proposition \ref{1.7prop}(ii).
\eit
{}From the same arguments as in Theorem \ref{fs}, we get a contradiction (for some nonzero primitive ideal $Q_i$, $x-\lambda+Q_i$ is, on the one hand, invertible in $R/Q_i$, and on the other hand, belongs to a maximal left ideal of $R/Q_i$).
Therefore, $R$ is algebraic over $k$.

\end{proof}

\subsection{Boundedness}
Let $R$ be an affine almost PI algebra.
For each ideal $0\neq I \triangleleft R$ there exists $n=n(I)$, which depends on $I$, such that $J(R/I)^{n(I)}=\bar{0}$.
This follows from the Razmyslov-Kemer-Braun theorem \cite{belov} which says that every affine PI algebra has a nilpotent Jacobson radical.
This leads to the following definition.

\begin{defin}\label{nilpotent-bounded}
A ring $R$ is \emph{bounded}, if there exists some $d$ such that for every $0\neq I\triangleleft R$, $J(R/I)^d=\bar{0}$.

Any such $d$ will be called a \emph{nilpotency-bound} for $R$.
\end{defin}
An affine almost PI algebra may or may not be bounded (when there exists no such a $d$).
Actually, boundedness of an affine almost PI algebra is a very strong condition, as the following theorem shows.
A theorem of Passman and Temple or of Farkas, helps us to prove it.
\begin{thm}\label{bounded}
If an affine almost PI algebra $R$ is bounded, then it is non-PI just infinite.
\end{thm}
\begin{proof}
Take $0\neq I\triangleleft R$. We must show that $\dim_k(R/I)<\infty$.

Let $d$ be a nilpotency-bound for $R$.
Theorem 3.3 in Passman and Temple \cite{PassTem} gives what we want, taking their $A$ to be our $R/I$ and their $d$ to be our bound $d$.

Theorem 3.3 in \cite{PassTem} says the following: "Let $A$ be an affine $k$ algebra satisfying a polynomial identity, and let $d$ be a fixed integer. If the nil radical of every finite dimensional homomorphic image of $A$ has nilpotence degree at most $d$, then $A$ itself is finite dimensional".

Indeed our $R/I$ is affine PI.
The nil radical of every homomorphic image of $R/I$ has nilpotence degree at most $d$.
This is because the nil radical is contained in the Jacobson radical, and the Jacobson radical has nilpotence degree at most $d$, as can easily be seen:
For $0\neq I\subseteq T\triangleleft R$ $(R/I)/(T/I)\cong R/T \Rightarrow \ J((R/I)/(T/I))\cong\ J(R/T) \Rightarrow \ J((R/I)/(T/I))^d\cong\ J(R/T)^d=\bar{0} \Rightarrow \ J((R/I)/(T/I))^d=\bar{0}$
\end{proof}

Theorem \ref{bounded} shows that if an affine almost PI algebra has a nonzero ideal with infinite codimension, then it cannot be bounded. This shows that boundedness is not a property that we should expect from an affine almost PI algebra to have.

\begin{cor}\label{cor.bounded}
If an affine almost PI algebra $R$ is bounded semiprimitive and the base field is uncountable, then $R$ is primitive.
\end{cor}
\begin{proof}
$R$ is non-PI just infinite from Theorem \ref{bounded}. Now use Theorem 2.2 in \cite{fs}
\end{proof}

\begin{cor}
If an affine almost PI algebra $R$ is bounded, then its center, $Z(R)$, is a field.

If, in addition, $R$ is not simple, then $\dim_kZ(R)<\infty$.
\end{cor}
\begin{proof}
Since $R$ is non-PI just infinite, we can use Proposition 2.3 in \cite{stably}.
\end{proof}

We bring now another proof how boundedness of $R$ implies $Z(R)$ is a field, but with no information about $\dim_kZ(R)$.
We include it here, because it shows that for every nonzero prime $P$, $Z(R/P)$ also a field. 
However, this is already known, since $R/P$ is simple (from Theorem \ref{bounded}, $R$ is just infinite, so $R/P$ is prime finite dimensional, hence simple).

\begin{thm}\label{bounded2}
If an affine almost PI algebra $R$ is bounded, then $Z(R)$ is a field.
Moreover, for every nonzero prime $P$, $Z(R/P)$ is also a field.
\end{thm}

\begin{proof}
The following is an argument applicable to any prime algebra $S$ whose homomorphic images (except, maybe, $S/0$) has a nilpotent Jacobson radical and there exists a bound of nilpotency.
We will prove that for such algebra $S$, its center is a field.

Since $S$ is prime, each $0\neq z\in Z(S)$ is regular.
Take $0\neq z\in Z(S)$.
It is enough to prove that $z$ has an inverse in $S$.
Indeed, assume that $y\in S$ s.t. $zy=yz=1$.
For every $r\in S$, $z(yr-ry)=z(yr)-z(ry)=(zy)r-r(zy)=1r-r1=0$, so by regularity of $z$ we have $yr-ry=0$, hence $y\in Z(S)$.
\bit
\item Step 1: For each $n\geq1$ we have $z+Sz^n\in J(S/Sz^n)$:
Fix $n$.
Remember $J(S/Sz^n)=J(\bar{S})=\{\bar{w}\in \bar{S} | \ \text{for every}\ \bar{r}\in \bar{S}: \ \bar{1}+\bar{r}\bar{w} \ \text{is left invertible} \}$.
Let us see why for every $\bar{r} \in \bar{S}$: $\bar{1}+\bar{r}\bar{z}$ is left invertible.
Take $\bar{r}\in \bar{S}$.
Define $y_r:=1-rz+r^2z^2- \ldots \pm r^{n-1}z^{n-1}$. 
$\bar{y_r}(\bar{1}+\bar{r}\bar{z})=\overline{y_r(1+rz)}=\overline{(1-rz+r^2z^2- \ldots \pm r^{n-1}z^{n-1})(1+rz)}=\overline{1\pm r^nz^n}=\bar{1}$,
so $\bar{y_r}$ is the left inverse of $\bar{1}+\bar{r}\bar{z}$.
(Notice that for $n=1$: $\bar{1}+\bar{r}\bar{z}=\overline{1+rz}=\bar{1}$, so it is trivial that it has a left inverse,
namely itself. Actually , we will need the above only in case $n>1$).

\item Step 2: Assume $d\geq1$ is a bound of nilpotency for $S$ and take $m>d$.
From step 1, in particular for that $m$: $z+Sz^m\in J(S/Sz^m)$.
Then, $z^d+Sz^m=(z+Sz^m)^d\in J(S/Sz^m)^d=\bar{0}=0+Sz^m$, so $z^d\in Sz^m$.
Hence, $z^d=rz^m$ for some $r\in S$, and we have $(1-rz^{m-d})z^d=z^d-rz^m=0$.
Finally, since $z$ is regular $1-rz^{m-d}=0$, so $rz^{m-d}=1$ and $z$ is invertible in $S$.
\eit
For the second argument: Take a nonzero prime ideal $P$. 
$R/P$ has the property that each of its homomorphic images has nilpotent Jacobson radical and there exists a bound of nilpotency (the bound for $R$ holds for $R/P$), so we can use the above argument to show that $Z(R/P)$ is a field.
\end{proof}

\begin{cor}
If an affine almost PI algebra $R$ is bounded, then each nonzero prime ideal is maximal.
\end{cor}
\begin{proof}
Actually, this follows immediately from Theorem \ref{bounded}, which says that $R$ is just infinite.
Indeed, if $P$ is a nonzero prime ideal, then $R/P$ is prime finite dimensional, hence simple.

Let us present another argument:
Let $P$ be any nonzero prime ideal in $R$. From Theorem \ref{bounded2}, $Z(R/P)$ is a field.
Since $R/P$ is prime PI, it follows from Rowen's Theorem that $R/P$ is simple, so $P$ is maximal.
\end{proof}

A somewhat similar idea is that of degree-boundedness:
\begin{defin}
An almost PI algebra (not necessarily affine) is \emph{degree-bounded}, if there exists some $d$ such that for every $0\neq I\triangleleft R$, $R/I$ satisfies a polynomial identity of degree $d$.
Any such $d$ will be called a \emph{degree-bound} for $R$
\end{defin}

However, degree-boundedness seems too strong to require, as can, for example, be seen from the following proposition.
\begin{prop}
Let $R$ be an almost PI algebra. If $R$ is degree-bounded and semiprimitive, then it is primitive.
\end{prop}

\begin{proof}
Otherwise, 0 is not primitive, so $J(R)=\cap Q_\lambda$ where each $Q_\lambda$ is a nonzero primitive ideal.
But $J(R)=0$, hence $R\hookrightarrow \prod R/Q_\lambda$, a contradiction (reason: if $d$ is a degree-bound for $R$, then $\prod R/Q_\lambda$ satisfies a polynomial identity of degree $d$, hence $R$ satisfies a polynomial identity of degree $d$, which is impossible).
\end{proof}

\subsection{Stably almost PI}
This subsection is based on \cite{stably}. Notice that affinity is not always needed.
Luckily, some (but not all) of the theorems in \cite{stably} which deals with a non-PI just infinite algebra are still true if we take an almost PI algebra instead.

\begin{defin}
Let $R$ be an almost PI algebra over a field $k$. $R$ is \emph{stably almost PI} if the $K$-algebra $R\otimes_k K$ is an almost PI algebra, for every field extension $K/k$.
\end{defin}
From now on, let $C(R)$ denote the extended center of $R$ and $Z(R)$ denote the center of $R$.
\begin{thm}\label{1.41}
Let $R$ be an almost PI algebra over an algebraically closed field $k$ with $Z(R)=k$. Then $R$ is stably almost PI $\Leftrightarrow$ $R$ is centrally closed (this means $C(R)=Z(R)$).
\end{thm}
\begin{proof}
Reading carefully the proof of Theorem 3.4 in \cite{stably} reveals that their proof also works here. Of course minor changes must be done in their proof: When $R$ is simple it is centrally closed. $R$ is stably almost PI since $R\otimes_k K$ is simple for every field extension $K/k$, hence almost PI.

In the notations of \cite{stably} (the only difference in notations is that our $R$ is their $A$):

($\Leftarrow$) When $R$ is centrally closed $J$ is still a nonzero ideal of $R$. $R/J$ is PI (since $R$ is almost PI) and $K$ is PI, therefore, from \cite[Theorem 2.4]{gol}, $(R/J)\otimes_k K$ is PI. Since $(R\otimes_k K)/I \stackrel{onto}{\leftarrow} (R\otimes_k K)/(J\otimes_k K)\cong (R/J)\otimes_k K$, we see that $(R\otimes_k K)/I$ is PI (as a homomorphic image of a PI ring).

($\Rightarrow$) Supposing $R$ is stably almost PI but not centrally closed, all the arguments in \cite{stably} are valid here untill the last three lines: $R$ embeds in $B/(a+tb)$ ($B:=R\otimes_k k(t)$). $R$ is stably almost PI so $B/(a+tb)$ is PI. Thus, $R$ is PI, a contradiction.
\end{proof}

As in \cite{stably}, one gets the following propositions:
\begin{prop}
Let $R$ be an almost PI $k$-algebra. Then $R$ is stably almost PI $\Leftrightarrow$ the $C(R)$-algebra $R\otimes_k C(R)$ is almost PI.
\end{prop}

\begin{proof}
See Proposition 3.5 in \cite{stably}.
\end{proof}

\begin{prop}
Let $R$ be a primitive almost PI $k$-algebra, where $k$ is an algebraically closed field. If $R$ satisfies the Nullstellensatz, then $R$ is stably almost PI.
\end{prop}

\begin{proof}
See Proposition 4.3 in \cite{stably}.
\end{proof}

Next we consider countably generated almost PI algebras.

The following proposition is immediate from a theorem of Rowen and Small.

\begin{prop}\label{1.7.32}
Let $R$ be a countably generated almost PI $k$-algebra. If $R$ is Goldie and $k$ is uncountable, then $C(R)$ is algebraic over $k$.
\end{prop}

\begin{proof}
Use Theorem 1 in \cite{goldie}, which says the following: "Suppose $R$ is a prime Goldie algebra over a field $k$. If $\dim_k(R)<|k|$, and the nonzero prime ideals of $R$ have a separating set of some cardinality $\gamma<|k|$, then $C(R)$ is algebraic over $k$". Of course our $R$ is prime (Theorem \ref{prime}), $\dim_k(R)=\aleph_0<|k|$, and the nonzero (prime) ideals of $R$ have a separating set of cardinality $\gamma=\aleph_0<|k|$ (Proposition \ref{ccsubset}).
\end{proof}

We can remove the assumption that $R$ is Goldie and still conclude that $C(R)$ is algebraic over $k$. Actually more can be said:

\begin{prop}\label{extended}
Let $R$ be a countably generated almost PI $k$-algebra. Then $C(R)$ is a countable dimensional field extension of $k$. Therefore, if $k$ is uncountable, then $C(R)$ is algebraic over $k$.
\end{prop}

\begin{proof}
From Proposition 6.2 and Remark 6.3 in \cite{stably}, $C(R)$ is a countable dimensional field extension of $k$. 
If $k$ is uncountable, Amitsur's trick yields that $C(R)$ is algebraic over $k$ (remember that $C(R)$ is a field, hence for each $x\in C(R)-k$, $\{x-\lambda1 \}_{\lambda\in k}$ are all invertible, etc.).
\end{proof}

\begin{cor}\label{44}
Let $R$ be a countably generated almost PI $k$-algebra and let $k$ be an uncountable algebraically closed field. Then $R$ is stably almost PI.
\end{cor}

\begin{proof}
From Proposition \ref{extended}, $C(R)$ is algebraic over $k$. But $k$ is algebraically closed, hence $C(R)=k$. Therefore $C(R)=Z(R)=k$, so $R$ is centrally closed. Theorem \ref{1.41} shows that $R$ is stably almost PI.
\end{proof}

The uncountable hypothesis in Corollary \ref{44} is necessary, since there exists an affine almost PI algebra over a countable algebraically closed field $k$ which is not stably almost PI; just take the algebra $A$ from the example in \cite{stably}. Indeed, the algebra $A$ is non-PI just infinite affine over a countable algebraically closed field $k$ and it is not centrally closed. Hence $A$ is not stably almost PI by Theorem \ref{1.41}.

Notice that Proposition \ref{1.7.32}, Proposition \ref{extended} and Corollary \ref{44} can be generalised to "the usual dimension-cardinality hypothesis":

\begin{prop}[generalised Proposition 1.44]
Let $R$ be an almost PI $k$-algebra with $\dim_k(R)<|k|$. If $R$ is Goldie, then $C(R)$ is algebraic over $k$.
\end{prop}

\begin{proof}
Use Theorem 1 in \cite{goldie}. Observe that indeed, the nonzero (prime) ideals of $R$ have a separating set of cardinality  $\gamma=\aleph_0\leq \dim_k(R)<|k|$ (Proposition \ref{ccsubset}).
\end{proof}

Again, we can remove the assumption that $R$ is Goldie and still conclude that $C(R)$ is algebraic over $k$:

\begin{prop}[generalised Proposition 1.45]\label{48}
Let $R$ be an almost PI $k$-algebra. If $\dim_k(R)<|k|$, then $C(R)$ is algebraic over $k$.
\end{prop}

\begin{proof} 
Use Proposition 6.2 in \cite{stably} to get $\dim_kC(R)\leq \dim_k(R)$. If $\dim_k(R)<|k|$, then $\dim_kC(R)\leq \dim_k(R)<|k|$ and Amitsur's trick implies that $C(R)$ is algebraic over $k$.
\end{proof}

\begin{cor}[generalised Corollary 1.46]
Let $R$ be an almost PI $k$-algebra. If $\dim_k(R)<|k|$ and $k$ is an algebraically closed field, then $R$ is stably almost PI.
\end{cor}

\begin{proof}
From Proposition \ref{48}, $C(R)$ is algebraic over $k$. But $k$ is algebraically closed, hence $C(R)=k$. Therefore $C(R)=Z(R)=k$, so $R$ is centrally closed. Theorem \ref{1.41} shows that $R$ is stably almost PI.
\end{proof}

Finally, Proposition 2.3 in \cite{stably} says the following: "Let $A$ be a non-PI just infinite $k$-algebra that is not simple. Then $Z(A)$ is a finite dimensional field extension of $k$". 

In the proof there, after knowing that $Z(A)$ is a field, non-simplicity of $A$ guarantees that $\dim_k(Z(A))<\infty$.
 
For an almost PI algebra we have a similar result, but the proof is different. Also notice that for an almost PI algebra we assume, in addition, that $\dim_k(R)<|k|$.

\begin{prop}\label{centerfdfe}
Let $R$ be an almost PI $k$-algebra with $\dim_k(R)<|k|$. Then $Z(R)$ is a field extension of $k$.

If $R$ has a nonzero proper ideal which has finite codimension (for example: if $R$ is affine and not simple), then $Z(R)$ is a finite dimensional field extension of $k$.
\end{prop}

\begin{proof}
Since $R$ is an almost PI $k$-algebra with $\dim_k(R)<|k|$, $C(R)$ is algebraic over $k$ (Proposition \ref{48}); hence $Z(R)\subseteq C(R)$ is algebraic over $k$. Primality of $R$ (Theorem \ref{prime}) implies $Z(R)$ is a domain. Therefore, $Z(R)$ is a field. 

Let $I$ be a nonzero proper ideal of $R$ with $\dim_k(R/I)< \infty$. Since $Z(R)$ is a field, $I\cap Z(R)=0$. Therefore, the natural projection $R\longrightarrow R/I$ induces an injection $Z(R)\hookrightarrow R/I$, so $\dim_kZ(R)\leq \dim_k(R/I)< \infty$.

In the special case where $R$ is affine and not simple: From Zorn's Lemma, there exists a maximal ideal $Q$. $Q$ is nonzero (otherwise, $R$ would be simple). As an affine almost PI algebra $\dim_k(R/Q)< \infty$ (remember Proposition \ref{1.7prop}(ii)).
\end{proof}

\begin{cor}
Let $R$ be an almost PI $k$-algebra with $\dim_k(R)<|k|$ and $k$ is algebraically closed. 
If $R$ has a nonzero proper ideal which has finite codimension (for example: if $R$ is affine and not simple), then $R$ is central ($Z(R)=k$).
\end{cor}

\begin{proof}
Immediate from the above proposition.
\end{proof}

\section{An affine prime Goldie algebra with GK dimension <3 over an uncountable field}
Let $R$ be an affine prime Goldie $k$-algebra with GK dimension $<3$, where $k$ is an uncountable field. We want to know when such an algebra is PI or primitive.

If such an $R$ is algebraic, then it is simple.
This is because any algebraic prime Goldie algebra is simple; Indeed, in a prime Goldie algebra each nonzero two sided ideal contains a regular element. This regular element is algebraic, hence invertible (actually, we have already seen that explanation in Corollary \ref{1.cor.simple}(ii)).

If such an $R$ is von Neumann regular, then it is simple.
Indeed, a Goldie von Neumann regular algebra is Artinian, see \cite[Corollary 2.16]{goodvnr}, and a prime Artinian algebra is simple, see \cite[2.3.9]{ro1}.

In the following subsections we will bring some partial answers to Small's question, namely, that a non-PI such $R$ is primitive when:
\bit
\item [(1)] There exists a non-algebraic element $x \in R$, such that $x+P \in R/P$ is algebraic over $k$, for every prime ideal $P$ with co-GK 1, and $k$ is algebraically closed.

\item [(2)] There exists a non-algebraic element $x \in R$, such that $x+P \in R/P$ is algebraic over $k$, for every \emph{completely} prime ideal $P$ with co-GK 1, 
and $x+P \in R/P$ is nilpotent, for every prime ideal $P$ with co-GK 1, which is not completely prime.

\item [(3)] There exists a non-algebraic element $x \in R$, such that $x+P \in R/P$ is algebraic over $k$, for every prime ideal $P$ with co-GK 1, and also for every such prime $P$, one of the following two conditions is satisfied:
\bit
\item the polynomial ring $Z(R/P)[t]$ has the property (*) \ref{property}.
\item the center of the classical (Goldie) ring of quotients of $R/P$, $Z(Q(R/P))$, is isomorphic to $k(X)$.
\eit
\eit

First we give two easy lemmas.
\begin{lem}\label{2.semiprimitive}
Let $R$ be a prime Goldie $k$-algebra with $dim_k(R)<|k|$. Then $R$ is semiprimitive.
\end{lem}

\begin{proof}
$\dim_k(R)<|k|$ implies $J(R)$ is algebraic over $k$ (use Amitsur's trick).
Prime Goldie implies that each nonzero two sided ideal contains a regular element.

Hence, if $J(R) \neq 0$, then $J(R)$ contains a regular element, $x$. That $x$ is also algebraic over $k$. Therefore, $x$ is invertible, which is impossible.
\end{proof} 

\begin{lem}\label{lemma}
Let $R$ be an affine prime Goldie (non-simple) non-PI $k$-algebra with GK dimension $<3$, where $k$ is an uncountable field.
Then 
\bit
\item [(a)] $R$ is an almost PI algebra.
\item [(b)] $R$ has only countably many height 1 primes.
\item [(c)] a nonzero prime ideal $P$ can be one of the following three kinds:
\bit
\item [(1)] First kind: $\Ht(P)=2$, $\GKdim(R/P)=0$: P is necessarily maximal.
\item [(2)] Second kind: $\Ht(P)=1$, $\GKdim(R/P)=0$: P is necessarily maximal.
\item [(3)] Third kind: $\Ht(P)=1$, $\GKdim(R/P)=1$: P is necessarily not primitive.
\eit
\eit
\end{lem}

\begin{proof}
\bit
\item [(a)] Take a nonzero ideal $I$ of $R$. Since $R$ is prime Goldie $I$ contains a regular element. Therefore, from \cite[Proposition 3.15]{kra}, $\GKdim(R/I)\leq \GKdim(R)-1$. Hence, from Bergman's gap theorem \cite[Theorem 2.5]{kra}, $\GKdim(R/I)\in\{0,1\}$, so $R/I$ is PI (use \cite{ssw}).
\item [(b)] Use (a) and Proposition \ref{aa}.
\item [(c)] This follows from Corollary 3.16 in \cite{kra} which says: "Let $A$ be a $k$-algebra whose prime factor rings are right Goldie. If $P$ is a prime ideal of $A$, then $\GKdim(A)\geq \GKdim(A/P)+\Ht(P)$". 
Indeed we can use that corollary since; (i) from (a) $R$ is an almost PI algebra, hence prime (Theorem \ref{prime}) and we have assumed that $R$ is Goldie. (ii) a prime PI ring is Goldie (\cite[13.6.6(i)]{mr}).

Let $P$ be a prime of the first or of the second kind. Then $R/P$ as a prime finite dimensional algebra is simple, so $P$ is maximal.
Let $P$ be a prime of the third kind. $P$ is necessarily not primitive; otherwise, as a nonzero primitive ideal in an affine almost PI algebra, it must have finite codimension (see \ref{1.7prop}), a contradiction. 
\eit
\end{proof}

Notice that as a corollary to Lemma \ref{lemma}, one gets Lemma 2.8 of Bell's paper \cite{bell} which says the following: "Let $k$ be an uncountable field and let $A$ be a finitely generated prime Goldie $k$-algebra of GK dimension two that does not satisfy a polynomial identity. Then $A$ has at most countably many height one prime ideals $Q$ with the property that $A/Q$ is finite dimensional as a $k$ vector space". Indeed, from Lemma \ref{lemma}(b), $A$ has only countably many height 1 primes. Also notice that Lemma \ref{lemma}(b) tells us a little more; not only that there are only countably many height 1 primes with finite codimension, but also that there are only countably many height 1 primes with infinite codimension. This fact is crucial for the proof of Theorem \ref{algebraic images} and for the proof of Theorem \ref{nilpotent images}. Bell's lemma is enough for Theorem \ref{finitecoGK1}, since there we assume that there are only finitely many height 1 primes with infinite codimension.

Another corollary is the following, which will be used in the third section. 

\begin{cor}\label{2.cor}
Let $R$ be as in Lemma \ref{lemma}. If $R$ has a third kind prime, then there are uncountably many first kind primes.
\end{cor}

\begin{proof}
Otherwise, there are only countably many first kind primes.
Hence, for a third kind prime $P$ there are only countably many first kind primes which contain it, denote them $\{\tilde{Q_i}\}_{i \in I}$ where $|I|\leq \aleph_0$. Obviously, $\Spec(R/P)=\{\bar{0}\}\cup\{(\tilde{Q_i}/P)\}_{i \in I}$ with each $\tilde{Q_i}/P$ a maximal ideal. So $R/P$ has a countable separating set for its nonzero primitive ideals; just take any nonzero $x_i \in \tilde{Q_i}/P$. Actually, since $R/P$ is prime PI, we can choose nonzero $y_i \in \tilde{Q_i}/P \cap Z(R/P)$. But then, a theorem of Rowen \cite[Proposition 3.2]{rosem} implies that $R/P$ is primitive, a contradiction ($P$, as a third kind prime, is necessarily not primitive).
\end{proof}

\subsection{There are only finitely many prime ideals with co-GK 1}
Let $R$ be an affine prime Goldie non-simple non-PI $k$-algebra with GK dimension $<3$, where $k$ is an uncountable field.
One can ask what happens in the special case when $R$ has only finitely many prime ideals $P$ with $\GKdim(R/P)=1$ (since such primes are of height 1, the above Lemma \ref{lemma}(b) shows that there are only countably many of them).
The answer is: $R$ must be primitive.
This answer is already known and brought in Bell's paper \cite[page 832]{bell}. However, we bring its proof now.

\begin{thm}\label{finitecoGK1}
Let $R$ be an affine prime Goldie non-simple non-PI $k$-algebra with GK dimension $<3$, where $k$ is an uncountable field.
Asuume that $R$ has only finitely many prime ideals $P$ with $\GKdim(R/P)=1$.
Then $R$ is primitive.
\end{thm}

\begin{proof}
Lemma \ref{lemma} tells us that $R$ can have three kinds of nonzero primes; we shall call them primes of the first, second and third kind, as in Lemma \ref{lemma}.
So our assumption here is that $R$ has only finitely many third kind primes.

First option: There are no third kind primes.
Since there are no third kind primes, it follows that there are no first kind primes (a first kind prime must contain a third kind prime, but there are no third kind primes).
So, $\Spec(R)=\{0, \text{second kind primes}\}$.
Now we can apply \cite[Corollary 1]{prime} which says that if in an affine infinite dimensional algebra each nonzero prime ideal has finite codimension, then the algebra is just infinite, and get that $R$ is just infinite.
{}From Lemma \ref{2.semiprimitive} our $R$ is semiprimitive, so from Farkas and Small's theorem \cite[Theorem 2.2]{fs} (which says that an affine semiprimitive just infinite algebra over an uncountable field is PI or primitive), it follows that $R$ is primitive.

Second option: There is at least one third kind prime.
By assumption, there are only finitely many third kind primes, denote them  $\{P_1,\ldots,P_d\}$. 
Denote the first kind primes by $\tilde{T}=\{\tilde{Q_\varphi}\}_{\varphi \in \Phi}$ (From Corollary \ref{2.cor} $|\Phi|>\aleph_0$).
Denote the second kind primes by $T=\{Q_i\}_{i \in M}$. From Lemma \ref{lemma}(b) $|M|\leq \aleph_0$.

(a) $|M|<\aleph_0$: Then it is immediate that $R$ is primitive;

Otherwise, $0=J(R)=(\cap \tilde{Q_\varphi})\cap(\cap Q_i)\supseteq (P_1 \cap \cdots \cap P_d)\cap(\cap Q_i)\neq 0$, where $(P_1 \cap \cdots \cap P_d)\cap(\cap Q_i)\neq 0$ since it is a finite intersection of nonzero ideals in a prime ring, a contradiction. So $0$ is a primitive ideal of $R$.

(b) $|M|=\aleph_0$:
The nonzero primitives are exactly the primes of the second and of the first kind: $T\cup\tilde{T}=\{Q_i\}_{i \in M}\cup\{\tilde{Q_\varphi}\}_{\varphi \in \Phi}$.

$R$ is prime so $P_1\cap\cdots\cap P_d\neq0$. As a nonzero two sided ideal of a prime Goldie ring, $P_1\cap\cdots\cap P_d$ contains a regular element, call it $y$. This $y$ is not algebraic over $k$ (Otherwise, $y$ is regular and algebraic, hence invertible. But then $P_1\cap\cdots\cap P_d=R$, a contradiction).

Define two countable subsets of $k$:

$A:=\{\alpha \in k | \ y-\alpha \text{ is left invertible in R}\}$: Indeed, from Amitsur's trick $A$ is countable ($y$ is non-algebraic).

$B:=\{\beta \in k | \ \text{there exists}\ i \in M \text{ such that} \ y-\beta+Q_i \in R/Q_i \\
\text{ is not invertible in } R/Q_i\}$: Indeed, $\forall i \in M$ $R/Q_i\hookrightarrow \End_k(R/Q_i)\cong \M_{\dim_k(R/Q_i)}(k)$.
Therefore, there exists only finitely many $\beta \in k$ (=the roots of the characteristic polynomial of $y+Q_i \in R/O_i \hookrightarrow \End_k(R/Q_i)\cong \M_{\dim_k(R/Q_i)}(k)$) such that $y-\beta+Q_i=(y+Q_i)-(\beta+Q_i)$ is not invertible in $R/Q_i$, denote them $B(i)=\{\beta(i)_1,\ldots,\beta(i)_{m_i}\}$. Then $B=\cup_{i \in M}B(i)$, which is countable.

Next, define $E:=k-(A\cup B)$. $E$ is uncountable.

Take any nonzero $\lambda \in E$.

$y-\lambda$ is not left invertible in $R$ (since $\lambda \notin A$).
Hence $R(y-\lambda)\subsetneq R$, so there exists a maximal left ideal $L$ such that $R(y-\lambda)\subseteq L$.
$\ann_R(R/L)$ is a primitive ideal of $R$.

There are two options for $\ann_R(R/L)$:
\bit
\item [(I)] $\ann_R(R/L)=0$. Then $R$ is primitive.
\item [(II)] $\ann_R(R/L) \neq 0$. 
Hence $\ann_R(R/L) \in T\cup \tilde{T}$.

We shall now see that $\ann_R(R/L) \in T\cup \tilde{T}$ is impossible:
\bit
\item [(1)] $\ann_R(R/L) \in T$:
So $\ann_R(R/L)=Q_i$ for some $i \in M$.

$Q_i=\ann_R(R/L)=\core(L)\subseteq L$.
We have an impossible situation:
On the one hand, $y-\lambda+Q_i$ is invertible in $R/Q_i$ (since $\lambda \notin B$).
On the other hand,  $y-\lambda+Q_i\in L/Q_i\subseteq R/Q_i$, where $L/Q_i$ is a maximal left ideal of $R/Q_i$.
So $\ann_R(R/L) \in T$ cannot be true.
\item [(2)] $\ann_R(R/L) \in \tilde{T}$:
So $\ann_R(R/L)=\tilde{Q_\varphi}$ for some $\varphi \in \Phi$.
There exists $1 \leq j \leq d$ such that $P_j \subsetneq \tilde{Q_\varphi}=\ann_R(R/L)=\core(L)\subseteq L$.
We have an impossible situation:
On the one hand, $y-\lambda+P_j$ is invertible in $R/P_j$:  $y-\lambda+P_j=(y+P_j)+(-\lambda+P_j)=(0+P_j)+(-\lambda+P_j)=-\lambda+P_j$ (we have taken $y \in P_1\cap\cdots\cap P_d \subseteq P_j$).
On the other hand,  $y-\lambda+P_j\in L/P_j\subseteq R/P_j$, where $L/P_j$ is a maximal left ideal of $R/P_j$.
So $\ann_R(R/L) \in \tilde{T}$ cannot be true.
\eit
\eit
Hence only (I) is possible, so $R$ is primitive.
\end{proof}

As a special case of Theorem \ref{finitecoGK1} consider the following algebra. Remember that an ideal $P$ in a ring $R$ is \emph{completely prime}, if $R/P$ is a domain.

\begin{cor}
Let $R$ be an affine prime Goldie non-simple non-PI $k$-algebra with GK dimension $<3$, where $k$ is an uncountable algebraically closed field.
Assume that every prime ideal $P$ with $\GKdim(R/P)=1$ is completely prime.
Then $R$ is primitive.
\end{cor}

\begin{proof}
{}From Lemma \ref{lemma} $R$ has only countably many primes $P$ with $\GKdim(R/P)=1$.
There are two options:
\bit
\item[(i)] $R$ has only finitely many primes with co-GK 1.
\item[(ii)] $R$ has countably many (infinite) primes with co-GK 1.
But this option is impossible: An affine domain with GK dimension 1 over an algebraically closed field must be commutative (this claim is a consequence of a theorem of Small and Warfield \cite{sw} and a theorem of Tsen).
Hence, for every prime ideal $P$ with $\GKdim(R/P)=1$, $R/P$ is commutative.
From Remarks \ref{1.remarks} and Lemma \ref{1.ACCsemiprimes}, the intersection of those primes is zero. But then $R$ is commutative, a contradiction ($R$ is non-PI).
\eit
Therefore, only option (i) is possible. Now Theorem \ref{finitecoGK1} implies that $R$ is primitive.
\end{proof}

\subsection{There exists a non-algebraic element with algebraic images}
Continue to assume that $R$ is an affine prime Goldie non-simple non-PI $k$-algebra with GK dimension $<3$, where $k$ is an uncountable field.

In the following theorem \ref{algebraic images}, we show that $R$ is primitive, if, in addition:
\bit
\item $k$ is algebraically closed.
\item $R$ has a non-algebraic element $x \in R$ such that for every prime ideal $P$ with $\GKdim(R/P)=1$, $x+P \in R/P$ is algebraic over $k$.
\eit
If $R$ has only finitely many prime ideals $P$ with $GKdim(R/P)=1$, then $R$ is primitive even without the assumptions that $k$ is algebraically closed and that $R$ has a non-algebraic element $x$ such that for every prime ideal $P$ with $\GKdim(R/P)=1$, $x+P \in R/P$ is algebraic over $k$, as was proved in Theorem \ref{finitecoGK1}. So if that is the case, then we are done.

Therefore, our assumptions that $k$ is algebraically closed and that $R$ has such a special element are crucial when there are countably many (infinite) such primes. 

Notice:\bit
\item [(1)] When $R$ has only finitely many prime ideals $P$ with $\GKdim(R/P)=1$, then there exists a non-algebraic element $x \in R$ such that for every prime ideal $P$ with $\GKdim(R/P)=1$, $x+P \in R/P$ is algebraic over $k$:
(If $R$ has no such primes, then $R$ is just infinite, this was explained in the beginning of the proof of Theorem \ref{finitecoGK1}, in First option, and then $R$ is primitive).
Denote the primes with co-GK 1 by $\{P_1,\ldots,P_d\}$. Since $R$ is prime Goldie, there exists a regular element $x \in P_1 \cap\cdots\cap P_d$ which is necessarily non-algebraic. Obviously, for every $1\leq j \leq d$, $x+P_j=0+P_j$ is algebraic over $k$.
 
But if $R$ has infinitely many prime ideals $P$ with $GKdim(R/P)=1$ 
(necessarily countably many such primes, see Lemma \ref{lemma}(b)), 
then it seems to be yet unknown if $R$ must have such a special non-algebraic element $x$.

\item [(2)] In general, there exist non-algebraic algebras having a non-algebraic element whose image in every nonzero prime is algebraic. For example: $k[X]$ or any non-algebraic just infinite algebra. Indeed, if $x$ is a non-algebraic element in a just infinite algebra, then for every nonzero two sided ideal $I$ (prime or not), $x+I \in R/I$ is algebraic over $k$. ($R/I$ is finite dimensional, therefore every element in $R/I$ is algebraic, in particular $x+I$).

\item [(3)] We could have imposed a more restrictive assumption than the assumption that $R$ has a non-algebraic element $x$ such that for every prime ideal $P$ with $\GKdim(R/P)=1$, $x+P \in R/P$ is algebraic over $k$, namely that $R$ is an \textit{almost algebraic} algebra.

By definition, an almost algebraic algebra is an algebra $R$ over a field $k$ which is non-algebraic over $k$, but $R/I$ is algebraic over $k$, for every nonzero two sided ideal $I$ of $R$. Examples of an almost algebraic algebra are $k[X]$ or any non-algebraic just infinite algebra.

An almost algebraic algebra $R$ satisfies the above assumption: $R$ is non-algebraic over $k$; hence there exists a non-algebraic element $x \in R$. For every nonzero two sided ideal $I$ of $R$, $R/I$ is algebraic over $k$, so, in particular, $x+I \in R/I$ is algebraic.

(One can easily show that an almost algebraic algebra $R$ is prime: Let $x \in R$ be non-algebraic. If $R$ is not prime, then there exist $A,B$ nonzero ideals such that $AB=0$. But, $R/A$ is algebraic, so $f(x+A)=0+A$ for some $f(t) \in k[t]$, implying $f(x) \in A$. Similarly, $g(x) \in B$ for some $g(t) \in k[t]$. Therefore, $(fg)(x)=f(x)g(x) \in AB=0$, a contradiction to our choice of a non-algebraic $x$).

However, an affine almost PI algebra $R$ which is also almost algebraic is just infinite; Indeed, for every nonzero two sided ideal $I$ of $R$, $R/I$ is an algebraic affine PI algebra, hence finite dimensional (Kurosch problem for PI, see \cite[13.8.9(2)]{mr}).
\eit

\begin{thm}\label{algebraic images}
Let $k$ be an uncountable algebraically closed field and let $R$ be an affine prime Goldie non-simple non-PI $k$-algebra with GK dimension $<3$. 
If $R$ has a non-algebraic element $x$ such that $x+P \in R/P$ is algebraic over $k$, for every prime ideal $P$ with $GKdim(R/P)=1$, then $R$ is primitive.
\end{thm}

\begin{proof}
If there are only finitely many primes with co-GK 1, then we are done.

So assume that there are infinitely many such primes (necessarily countably many).
Denote them by $\{P_i\}_{i \in \N}$ (third kind primes). 
Denote the first kind primes by $\tilde{T}=\{\tilde{Q_\varphi}\}_{\varphi \in \Phi}$.
Denote the second kind primes by $T=\{Q_i\}_{i \in M}$. From Lemma \ref{lemma}(b) $|M|\leq \aleph_0$.
Define three countable subsets of $k$.

$A:=\{\alpha \in k | \ x-\alpha \text{ is left invertible in R}\}$: Indeed, from Amitsur's trick $A$ is countable (By assumption, $x$ is non-algebraic).

$B:=\{\beta \in k | \ \text{there exists}\ i \in M \text{ such that} \ x-\beta+Q_i \in R/Q_i \\
\text{ is not invertible in } R/Q_i\}$: Indeed, $\forall i \in M$ $R/Q_i\hookrightarrow \End_k(R/Q_i)\cong \M_{\dim_k(R/Q_i)}(k)$.

Therefore, there exists only finitely many $\beta \in k$ (=the roots of the characteristic polynomial of $x+Q_i \in R/O_i \hookrightarrow \End_k(R/Q_i)\cong \M_{\dim_k(R/Q_i)}(k)$) such that $x-\beta+Q_i=(x+Q_i)-(\beta+Q_i)$ is not invertible in $R/Q_i$, denote them $B(i)=\{\beta(i)_1,\ldots,\beta(i)_{m_i}\}$. Then $B=\cup_{i \in M}B(i)$, which is countable.

$C:=\{\gamma \in k | \ \text{there exists}\ i \in \N \text{ such that} \ x-\gamma+P_i \in R/P_i \\
\text{ is not invertible in } R/P_i\}$:
Indeed, $\forall i \in \N $, $R/P_i$ is a prime PI algebra, hence can be embedded in a matrix ring over a field extension of $k$: $R/P_i\hookrightarrow \M_{n_i}(L_i)$, where $L_i$ is a field extension of $k$.

We have assumed that $x+P_i \in R/P_i$ is algebraic over $k$, so let $f_i(t) \in k[t]$ be such that $f_i(x+P_i)=0+P_i$.
Now think of $x+P_i \in R/P_i \stackrel{\rho_i}{\hookrightarrow} \M_{n_i}(L_i)$ as a matrix over $L_i$.
Obviously, $f_i(\rho_i(x+P_i))=\rho_i((f_i(x+P_i))=\rho_i(0+P_i)=0$. Hence, the minimal polynomial of $\rho_i(x+P_i)$, call it $m_i(t) \in L_i[t]$, divides $f_i(t)$. Remember we have assumed that $k$ is algebraically closed, so clearly $m_i(t) \in k[t]$. But then the characteristic polynomial of $\rho_i(x+P_i)$, call it $c_i(t) \in L_i[t]$ must also be over $k$, $c_i(t) \in k[t]$.
Let $C(i)=\{\gamma(i)_1,\ldots,\gamma(i)_{m_i}\}$ be the roots of $c_i(t)$ (of course, $C(i)\subset k$).
Therefore, for any $\gamma \in k-C(i)$, $x-\gamma+P_i=(x+P_i)-(\gamma+P_i)$ is invertible in $R/P_i$; Indeed,
$x+P_i$ algebraic over $k$ implies that $x-\gamma+P_i=(x+P_i)-(\gamma+P_i)$ is algebraic over $k$ (this trivial claim is obviously true for any $\gamma \in k$, but immediately we will need the particular choice of $\gamma \in k-C(i)$).
$x-\gamma+P_i=(x+P_i)-(\gamma+P_i)$ algebraic over $k$ implies that the characteristic polynomial of $\rho_i(x-\gamma+P_i)$ is over $k$, call it $c_{i,\gamma}(t)$ (of course, the argument is exactly as the argument for $c_i(t)$ to be over $k$).

So, let $c_{i,\gamma}(t)=(-1)^{n_i}t^{n_i}+b_{n_i-1}t^{n_i-1}+\ldots+b_1t+b_0$ with $b_{n_i-1},\ldots,b_1,b_0 \in k$.
{}From Cayley-Hamilton theorem, $c_{i,\gamma}(\rho_i(x-\gamma+P_i))=0$: 
(*) $(-1)^{n_i}(\rho_i(x-\gamma+P_i))^{n_i}+b_{n_i-1}(\rho_i(x-\gamma+P_i))^{n_i-1}+\ldots+b_1(\rho_i(x-\gamma+P_i))+b_0=0$.
Of course, $b_0=\det(\rho_i(x-\gamma+P_i))$. 

Now, since $\gamma \in k$ is not a root of $c_i(t)$, $c_i(\gamma) \in k-0$.
So, $k-0 \ni c_i(\gamma)=\det(\rho_i(x+P_i)-\gamma I)
=\det(\rho_i(x+P_i)-\rho_i(\gamma+P_i))=(\det(\rho_i((x+P_i)-(\gamma+P_i))=\det(\rho_i(x-\gamma+P_i))=b_0$.

Therefore, (*) shows that $((-b_0)^{-1})((-1)^{n_i}(\rho_i(x-\gamma+P_i))^{n_i-1}+b_{n_i-1}(\rho_i(x-\gamma+P_i))^{n_i-2}+\ldots+b_1I)(\rho_i(x-\gamma+P_i))=I$, 
so the inverse of $\rho_i(x-\gamma+P_i)$ is

$((-b_0)^{-1})((-1)^{n_i}(\rho_i(x-\gamma+P_i))^{n_i-1}+b_{n_i-1}(\rho_i(x-\gamma+P_i))^{n_i-2}+\ldots+b_1I)$.

{}From this we can conclude that in $R/P_i$, $x-\gamma+P_i$ is invertible with inverse $((-b_0)^{-1})((-1)^{n_i}(x-\gamma+P_i)^{n_i-1}+b_{n_i-1}(x-\gamma+P_i)^{n_i-2}+\ldots+b_11) \in R/P_i$
(we use the fact that $\rho_i$ is injective).

Of course, for $\gamma \in C(i)$: $0=c_i(\gamma)=\det(\rho_i(x+P_i)-\gamma I)=\det(\rho_i(x+P_i)-\rho_i(\gamma+P_i))=(\det(\rho_i((x+P_i)-(\gamma+P_i))=\det(\rho_i(x-\gamma+P_i))$,
so $\rho_i(x-\gamma+P_i)$ is not invertible in $\M_{n_i}(L_i)$, which implies that $x-\gamma+P_i$ cannot be invertible in $R/P_i$.
Summarizing, $C=\cup_{i \in \N}C(i)$, which is countable.

Notice: The argument that $C$ is countable could also be used to show that $B$ is countable (but not vice versa).
We separated to two sets $B$ and $C$ (according to the two different options for a height 1 prime to have finite or infinite codimension, respectively), for two reasons; First, it seems better to give an easier argument where possible. 
Second, in Theorem \ref{nilpotent images} it will be clear that the argument that the set $C$ there (which is defined a little differently from the set $C$ here) is countable cannot necessarily be used to show that $B$ is countable.

Next, define $E:=k-(A\cup B \cup C)$. $E$ is uncountable.
Take any nonzero $\lambda \in E$.
$x-\lambda$ is not left invertible in $R$ (since $\lambda \notin A$).
Hence $R(x-\lambda)\subsetneq R$, so there exists a maximal left ideal $L$ such that $R(x-\lambda)\subseteq L$.
$\ann_R(R/L)$ is a primitive ideal of $R$.
There are two options for $\ann_R(R/L)$:
\bit
\item [(I)] $\ann_R(R/L)=0$. Then $R$ is primitive.
\item [(II)] $\ann_R(R/L) \neq 0$. 
Hence $\ann_R(R/L) \in T\cup \tilde{T}$.

We shall now see that $\ann_R(R/L) \in T\cup \tilde{T}$ is impossible:\bit
\item [(1)] $\ann_R(R/L) \in T$:
So $\ann_R(R/L)=Q_i$ for some $i \in M$.

$Q_i=\ann_R(R/L)=\core(L)\subseteq L$.
We have an impossible situation:
On the one hand, $x-\lambda+Q_i$ is invertible in $R/Q_i$ (since $\lambda \notin B$).
On the other hand,  $x-\lambda+Q_i\in L/Q_i\subseteq R/Q_i$, where $L/Q_i$ is a maximal left ideal of $R/Q_i$.
So $\ann_R(R/L) \in T$ cannot be true.
\item [(2)] $\ann_R(R/L) \in \tilde{T}$:
So $\ann_R(R/L)=\tilde{Q_\varphi}$ for some $\varphi \in \Phi$.
There exists $i \in \N$ such that $P_i \subsetneq \tilde{Q_\varphi}=\ann_R(R/L)=\core(L)\subseteq L$.
We have an impossible situation:
On the one hand, $x-\lambda+P_i$ is invertible in $R/P_i$ (since $\lambda \notin C$).  
On the other hand,  $x-\lambda+P_i\in L/P_i\subseteq R/P_i$, where $L/P_i$ is a maximal left ideal of $R/P_i$.
So $\ann_R(R/L) \in \tilde{T}$ cannot be true.
\eit
\eit
Hence only (I) is possible, so $R$ is primitive.
\end{proof}

\subsection{There exists a non-algebraic element with nilpotent images}
One may wish to take in Theorem \ref{algebraic images} an uncountable field which is not necessarily algebraically closed.
Going back to the proof of Theorem \ref{algebraic images} shows that the assumptions that $k$ is algebraically closed and $R$ has such a special non-algebraic element $x$, guarantees that for every height 1 prime $I$, the characteristic polynomial of $x+I \in R/I\hookrightarrow \M_{n_I}(L_I)$, $c_I(t)$, is over $k$.
(If $\GKdim(R/I)=1$, then $c_I(t)$ is apriori over $L_I$). 
This (and the fact that the characteristic polynomial of $x-\mu+I$ is over $k$, for any $\mu \in k$) guarantees that for every $\lambda \in k-D_I$, $x-\lambda +I$ is invertible in $R/I$ ($D_I$ is the set of roots of $c_I(t)$). 
Now, without assuming that $k$ is algebraically closed, we also want to find for every height 1 prime $I$, a finite (or countable) subset of $k$, call it $E_I$, such that for every $\lambda \in k-E_I$, $x-\lambda +I$ is invertible in $R/I$.
\bit
\item {}For the height 1 primes with finite codimension, this is immediate; just choose any non-algebraic element $x \in R$ and take $E_I$ to be the roots of the characteristic polynomial of $x+I \in R/I \hookrightarrow \End_k(R/I)\cong \M_{\dim_k(R/I)}(k)$.
\item {}For the height 1 primes with infinite codimension (co-GK 1), it is not yet clear whether algebraicity of an element $x+I \in R/I$ guarantees the existence of a finite (or countable) $E_I \subset k$, such that for every $\lambda \in k-E_I$, $x-\lambda +I$ is invertible in $R/I$.

Observe that now (when $k$ is not algebraically closed), even though the characteristic polynomial of $x+I \in R/I \hookrightarrow M_{n_I}(L_I)$ may not be over $k$, such a finite subset of $k$, $E_I$, may exist, as the following example shows:
Let $k=\R$, $R/I\hookrightarrow \M_2(\C(X))$ where $X$ is an indeterminate ($X$ commutes with $\C$) and $x+I=iI_2$.
$x+I$ is algebraic over $\R$ since it satisfies $t^2+1$, but its characteristic polynomial is $t^2-2it-1 \notin \R[t]$. 
However, one can take $E_I=\emptyset$; 
Indeed, for every $\lambda \in \R-E_I=\R$, $x-\lambda +I=(x+I)-(\lambda +I)= (iI_2)-(\lambda I_2) = (i-\lambda)I_2$, which is invertible in $R/I$, since its inverse is 
$(1/(i-\lambda))I_2= ((i+\lambda)/-(1+\lambda^2))I_2= (i+/-(1+\lambda^2))I_2 + (\lambda/-(1+\lambda^2))I_2 \in R/I$.
\eit

Again, as was mentioned in the discussion preceeding Theorem \ref{algebraic images}, if $R$ has only finitely many prime ideals $P$ with $GKdim(R/P)=1$, then $R$ is primitive without any further assumptions, as was proved in Theorem \ref{finitecoGK1}. So if that is the case, then we are done.

Therefore, the assumptions in the following Theorem \ref{nilpotent images} are crucial when there are countably many (infinite) such primes. 

\begin{thm}\label{nilpotent images}
Let $k$ be an uncountable field and let $R$ be an affine prime Goldie non-simple non-PI $k$-algebra with GK dimension $<3$. 
If $R$ has a non-algebraic element $x$ such that for every prime ideal $P$ with $GKdim(R/P)=1$:
\bit
\item [(i)] If $P$ is completely prime, $x+P \in R/P$ is algebraic over $k$
\item [(ii)] If $P$ is not completely prime, $x+P \in R/P$ is nilpotent,
\eit
then $R$ is primitive.
\end{thm}

\begin{proof}
If there are only finitely many primes with co-GK 1, then we are done.

So assume that there are infinitely many such primes (necessarily countably many).
Denote them by $\{P_i\}_{i \in \N}$ (third kind primes). Separate the third kind primes into two (disjoint) sets:
$\{P_i\}_{i \in N_1}$ the completely primes ($|N_1|\leq \aleph_0$), and $\{P_i\}_{i \in N_2}$ the ones that are not completely primes ($|N_2|\leq \aleph_0$). 
Denote the first kind primes by $\tilde{T}=\{\tilde{Q_\varphi}\}_{\varphi \in \Phi}$.
Denote the second kind primes by $T=\{Q_i\}_{i \in M}$. From Lemma \ref{lemma}(b) $|M|\leq \aleph_0$.
Define four countable subsets of $k$.

$A:=\{\alpha \in k | \ x-\alpha \text{ is left invertible in R}\}$: Indeed, from Amitsur's trick $A$ is countable (By assumption, $x$ is non-algebraic).

$B:=\{\beta \in k | \ \text{there exists}\ i \in M \text{ such that} \ x-\beta+Q_i \in R/Q_i \\
\text{ is not invertible in } R/Q_i\}$: Indeed, $\forall i \in M$ $R/Q_i\hookrightarrow \End_k(R/Q_i)\cong \M_{\dim_k(R/Q_i)}(k)$.

Therefore, there exists only finitely many $\beta \in k$ (=the roots of the characteristic polynomial of $x+Q_i \in R/O_i \hookrightarrow \End_k(R/Q_i)\cong \M_{\dim_k(R/Q_i)}(k)$) such that $x-\beta+Q_i=(x+Q_i)-(\beta+Q_i)$ is not invertible in $R/Q_i$, denote them $B(i)=\{\beta(i)_1,\ldots,\beta(i)_{m_i}\}$. Then $B=\cup_{i \in M}B(i)$, which is countable.

$C:=\{\gamma \in k | \ \text{there exists}\ i \in N_1 \text{ such that} \ x-\gamma+P_i \in R/P_i \text{ is not invertible in }\\R/P_i\}$: Fix $i \in N_1$. 
{}From (i), $x+P_i \in R/P_i$ is algebraic over $k$, hence for every $\gamma \in k$, $(x+P_i)-(\gamma+P_i)=x-\gamma+P_i$ is algebraic over $k$.
$R/P_i$ is a domain, hence for every $\gamma \in k$, except maybe one scalar $\gamma_i \in k$, $x-\gamma+P_i$ is regular (If $x+P_i$ is non-scalar, then for every $\gamma \in k$, $x-\gamma+P_i \neq \bar{0}$, hence regular. If there exists $\gamma_i \in k$ such that $x+P_i=\gamma_i+P_i$, then for every $\gamma \in k-\{\gamma_i\}$, $x-\gamma+P_i \neq \bar{0}$, hence regular). 
Therefore, for every $\gamma \in k$, except maybe one scalar $\gamma_i \in k$, $x-\gamma+P_i$ is invertible in $R/P_i$ (since an element which is algebraic and regular is invertible).
For every $i \in N_1$, let $C(i)=\{\gamma_i\}$ if $x+P_i=\gamma_i+P_i$, otherwise $C(i)=\emptyset$.
It is clear that $C=\cup_{i \in N_1}C(i)$, which is countable. 

$D:=\{\delta \in k | \ \text{there exists}\ i \in N_2 \text{ such that} \ x-\delta+P_i \in R/P_i \text{ is not invertible in }\\R/P_i\}$: Fix $i \in N_2$. 
$R/P_i$ is a prime PI algebra, hence can be embedded in a matrix ring over a field extension of $k$: $R/P_i\hookrightarrow \M_{n_i}(L_i)$, where $L_i$ is a field extension of $k$.

{}From (ii), $x+P_i \in R/P_i$ is nilpotent, so let $f_i(t)=t^{l_i}$ be such that $(x+P_i)^{l_i}=f_i(x+P_i)=0+P_i$

Now think of $x+P_i \in R/P_i \stackrel{\rho_i}{\hookrightarrow} \M_{n_i}(L_i)$ as a matrix over $L_i$.
Obviously, $f_i(\rho_i(x+P_i))=\rho_i((f_i(x+P_i))=\rho_i(0+P_i)=0$. Hence, the minimal polynomial of $\rho_i(x+P_i)$, call it $m_i(t) \in L_i[t]$, divides $f_i(t)$. 
But $f_i(t)=t^{l_i}$, so obviously $m_i(t)=t^{l'_i}$ with $l'_i \leq l_i$.
Then the characteristic polynomial of $\rho_i(x+P_i)$, call it $c_i(t) \in L_i[t]$ must be $c_i(t)=t^{n_i}$.
Take any nonzero $\delta \in k$. We claim that $x-\delta+P_i$ is invertible in $R/P_i$; 
Indeed, let $c_{i,\delta}(t) \in L_i[t]$ be the characteristic polynomial of $\rho_i(x-\delta+P_i)$.
$c_{i,\delta}(t)=\det(\rho_i(x-\delta+P_i)-tI)=\det(\rho_i((x+P_i)-(\delta+P_i))-tI)=
\det(\rho_i(x+P_i)-\rho_i(\delta+P_i)-tI)=\det(\rho_i(x+P_i)-\delta I-tI)=\det(\rho_i(x+P_i)-(\delta+t)I)=
c_i(\delta+t)=(\delta+t)^{n_i} \in k[t]$.
So, $c_{i,\delta}(t)=(\delta+t)^{n_i}=t^{n_i}+b_{n_i-1}t^{n_i-1}+\ldots+b_1t+b_0$ with $b_{n_i-1},\ldots,b_1,b_0 \in k$ (actually $b_j=C(n_i-j,n_i)\delta^{n_i-j}$. From Cayley-Hamilton theorem, 
$c_{i,\delta}(\rho_i(x-\delta+P_i))=0$. Then, it is clear that $\rho_i(x-\delta+P_i) \in \M_{n_i}(L_i)$ is invertible with inverse 
$((-b_0)^{-1})((\rho_i(x-\delta+P_i))^{n_i-1}+b_{n_i-1}(\rho_i(x-\delta+P_i))^{n_i-2}+\ldots+b_1I)$.

{}From this we can conclude that in $R/P_i$, $x-\delta+P_i$ is invertible with inverse $((-b_0)^{-1})((x-\delta+P_i)^{n_i-1}+b_{n_i-1}(x-\delta+P_i)^{n_i-2}+\ldots+b_11) \in R/P_i$
(we use the fact that $\rho_i$ is injective).

Define $D(i)=\{0\} \subset k$.
Of course, for $\delta \in D(i)$: $x-\delta+P_i=x+P_i$ is our nilpotent element, which is not invertible.
It is clear that $D=\cup_{i \in N_2}D(i)=\{0\}$, which is countable.

Next, define $E:=k-(A\cup B \cup C \cup D)$. $E$ is uncountable.
Take any $\lambda \in E$.
$x-\lambda$ is not left invertible in $R$ (since $\lambda \notin A$).
Hence $R(x-\lambda)\subsetneq R$, so there exists a maximal left ideal $L$ such that $R(x-\lambda)\subseteq L$.
$\ann_R(R/L)$ is a primitive ideal of $R$.
There are two options for $\ann_R(R/L)$:
\bit
\item [(I)] $\ann_R(R/L)=0$. Then $R$ is primitive.
\item [(II)] $\ann_R(R/L) \neq 0$.
Hence $\ann_R(R/L) \in T\cup \tilde{T}$.

We shall now see that $\ann_R(R/L) \in T\cup \tilde{T}$ is impossible:\bit
\item [(1)] $\ann_R(R/L) \in T$:
So $\ann_R(R/L)=Q_i$ for some $i \in M$.

$Q_i=\ann_R(R/L)=\core(L)\subseteq L$.
We have an impossible situation:
On the one hand, $x-\lambda+Q_i$ is invertible in $R/Q_i$ (since $\lambda \notin B$).
On the other hand,  $x-\lambda+Q_i\in L/Q_i\subseteq R/Q_i$, where $L/Q_i$ is a maximal left ideal of $R/Q_i$.
So $\ann_R(R/L) \in T$ cannot be true.

\item [(2)] $\ann_R(R/L) \in \tilde{T}$:
So $\ann_R(R/L)=\tilde{Q_\varphi}$ for some $\varphi \in \Phi$.
There exists $i \in \N=N_1 \cup N_2$ such that 
$P_i \subsetneq \tilde{Q_\varphi}=\ann_R(R/L)=\core(L)\subseteq L$.
We have an impossible situation:
On the one hand, $x-\lambda+P_i$ is invertible in $R/P_i$ (since $\lambda \notin (C \cup D)$).  
On the other hand,  $x-\lambda+P_i\in L/P_i\subseteq R/P_i$, where $L/P_i$ is a maximal left ideal of $R/P_i$.
So $\ann_R(R/L) \in \tilde{T}$ cannot be true.
\eit
\eit
Hence only (I) is possible, so $R$ is primitive.
\end{proof}

\subsection{Another partial answer to Small's question}
When $k$ is not algebraically closed, we have seen in Theorem \ref{nilpotent images} that the existence of a non-algebraic element $x \in R$ with algebraic images in completely primes with co-GK 1, and nilpotent images in the other primes with co-GK 1, guarantees that $R$ is primitive.
If one prefers the weaker assumption, that of Theorem \ref{algebraic images} when $k$ was algebraically closed,
(namely, the existence of a non-algebraic element $x \in R$ with algebraic images in primes with co-GK 1),
then one may add an additional condition in order to guarantee that $R$ is primitive.
The following definition will be used in the additional condition.
\begin{defin}\label{property}
Let $C$ be a commutative $k$-algebra ($k$ is any field).
{}For the polynomial ring $C[t]$, \emph{the property (*)}, is as follows:
If $f(t) \in C[t]$ and $g(t) \in C[t]$ have $f(t)g(t) \in k[t]$, then necessarily $f(t) \in k[t]$ and $g(t) \in k[t]$.
\end{defin}

For example: $C=k[X]$ has the property (*), since it is clear that if a product of two polynomials over $k[X]$ is a polynomial over $k$, then each of these polynomials must be over $k$.
A non-example is $C=k(X)$, since $1=X\cdot(1/X)$.

\begin{thm}\label{images}
Let $k$ be an uncountable field and let $R$ be an affine prime Goldie non-simple non-PI $k$-algebra with GK dimension $<3$. 
If $R$ has a non-algebraic element $x \in R$, such that $x+P \in R/P$ is algebraic over $k$, for every prime ideal $P$ with $\GKdim(R/P)=1$, and also for every such prime $P$, the polynomial ring $Z(R/P)[t]$ has the property (*),
then $R$ is primitive.
\end{thm}

\begin{proof}
The discussion in the beginning of the previous subsection shows that if the following condition is satisfied, then $R$ is primitive: For every prime ideal $P$ with $\GKdim(R/P)=1$, there exists a finite (or countable) subset of $k$, call it $E_P$, such that for every $\lambda \in k-E_P$, $x-\lambda+P$ is invertible in $R/P$.
Therefore, we only need to show that, in our $R$, for every prime ideal $P$ with $\GKdim(R/P)=1$, there exists such a finite set.

Let $P$ be a prime ideal with co-GK 1.
$R/P$ is a finite module over its center; this can be seen by either one of the following two arguments:
\bit
\item Small and Warfield have shown in \cite{sw} that an affine prime $k$-algebra with $\GKdim(R/P)=1$ is a finite module over its center.
\item Farina and Pendergrass-Rice have shown in \cite[Corollary 2]{prime} that a PI just infinite algebra is a finite module over its center (observe that $R/P$ is indeed a PI just infinite algebra).
\eit
Let $R/P\cong \Sigma_{i=1}^{n}Z(R/P)\bar{r_i}$, where $\bar{r_i} \in R/P$.
Hence, each element of $R/P$ can be thought of as a matrix over $Z(R/P)$. Reason: each $\bar{y} \in R/P$ can be thought of as an $Z(R/P)$-endomorphism of $R/P$, namely the left multiplication by $\bar{y}$.
The matrix which represents this left multiplication is, of course, $(a_{ij})$, where for each $1 \leq i \leq n$ $\bar{y}\bar{r_i}=\Sigma_{j=1}^{n}a_{ji}\bar{r_j}$. Indeed, $(a_{ij})[\bar{r_i}]=[\bar{y}\bar{r_i}]$).

(A remark: Although there may exist more then one way to think of an element of $R/P$ as a matrix over $Z(R/P)$ (since the $\bar{r_i}$'s are probably "dependent"), it seems that for what we would like to prove, this will not be a problem).

Now, in particular, $\bar{x}=x+P \in R/P\cong \Sigma_{i=1}^{n} Z(R/P)\bar{r_i}$ can be thought of as a matrix over $Z(R/P)$.
$Z(R/P)$ is a commutative domain, so it has a field of fractions, $F(Z(R/P))$.
$\bar{x}$ is algebraic over $k$, so let $f(t) \in k[t]$ be such that $f(\bar{x})=\bar{0}$. 
Let $\{A_j\}_{j \in J} \subseteq \M_n(Z(R/P))\subset \M_n(F(Z(R/P)))$ be all the possible options to view $\bar{x}$ as an element of $\M_n(Z(R/P))$ (we don't have any information about $|J|$, but we shall immediately see that this is not a problem).

Obviously, for each $A_j$: $f(A_j)=\bar{0}$ (since $f(\bar{x})=\bar{0}$ and if a matrix $B$ represents $\bar{x}$, then for any polynomial $h(t) \in Z(R/P)[t]$, $h(B)$ represents $h(\bar{x})$), so each $A_j$ is algebraic over $k$.
Denote by $c_j(t) \in F(Z(R/P))[t]$ the characteristic polynomial of $A_j \in \M_n(Z(R/P))\subset \M_n(F(Z(R/P)))$, and denote by $m_j(t) \in F(Z(R/P))[t]$ the minimal polynomial of $A_j$.

We claim that, for each $j \in J$, $c_j(t) \in k[t]$; Fix $j \in J$. In $F(Z(R/P))[t]$ we have: $f(t)=q(t)m_j(t)$ with $q(t) \in F(Z(R/P))[t]$, $m_j(t) \in F(Z(R/P))[t]$. From Gauss's lemma, $q(t) \in Z(R/P)[t]$, $m_j(t) \in Z(R/P)[t]$.
So, we have $f(t)=q(t)m_j(t)$ with $f(t) \in k[t]$ and $q(t),m_j(t) \in Z(R/P)[t]$.
But we have assumed that $Z(R/P)[t]$ has the property (*), therefore, $q(t),m_j(t) \in k[t]$.
Now, write $m_j(t)$ as a product of irreducible polynomials in $F(Z(R/P))$: $m_j(t)=irr_{1,j}(t)\cdots irr_{u_j,j}(t)$, where $irr_{1,j}(t),\ldots,irr_{u_j,j} \in F(Z(R/P))[t]$.
Use Gauss's lemma again to conclude that $irr_{1,j}(t),\ldots,irr_{u_j,j} \in Z(R/P)[t]$.
So, we have $m_j(t)=irr_{1,j}(t)\cdots irr_{u_j,j}(t)$ with $m_j(t) \in k[t]$ and $irr_{1,j}(t),\ldots,irr_{u_j,j} \in Z(R/P)[t]$.
But we have assumed that $Z(R/P)[t]$ has the property (*), therefore,  $irr_{1,j}(t),\ldots,irr_{u_j,j} \in k[t]$.
Finally, since each irreducible factor of $c_j(t)$ is an irreducible of $m_j(t)$, we have $c_j(t) \in k[t]$, as claimed.

Let $K$ be the algebraic closure of $k$, so $f(t)=(t-\alpha_1) \cdots (t-\alpha_d)$, where $d=deg(f)$ and $\alpha_1, \ldots ,\alpha_d \in K$ (it may happen that $\alpha_i=\alpha_j$ for $i \neq j$).
Define $\tilde{E_P}=:\{\alpha_1, \ldots ,\alpha_d \}$ ($\tilde{E_P} \subset K$).

Fix $j \in J$. Hence, $m_j(t)=(t-\alpha_{l_1}) \cdots (t-\alpha_{l_{s_j}})$, 
where $\{l_1,\ldots,l_{s_j}\} \subseteq \{1,\ldots,d\}$ and $s_j \leq d$.
So, $\alpha_{l_1},\ldots,\alpha_{l_{s_j}} \subseteq \tilde{E_P}$ are the only roots of $m_j(t)$.
Therefore, $\alpha_{l_1},\ldots,\alpha_{l_{s_j}} \subseteq \tilde{E_P}$ are also the only roots of $c_j(t)$.
This shows that the only possible options for the roots of any $c_j(t)$, $j \in J$, are among the elements of $\tilde{E_P}$.

Define $E_P:=\tilde{E_P}\cap k$.
For any $\lambda \in k-E_P$:
Fix $j \in J$. $c_j(\lambda)=det(\lambda I-A_j)\in k-0$ ($\in k$ since $c_j(t) \in k[t]$, and $\neq 0$ since $\lambda$ is not a root of $c_j(t)$).
Now, $\bar{x}-\lambda1 \in R/P$ is of course algebraic over $k$ (since $\bar{x} \in R/P$ is algebraic over $k$), so let $g(\bar{x}-\lambda1)=0$, with $g(t) \in k[t]$.
Let $\{B_l\}_{l \in L} \subseteq \M_n((Z(R/P))\subset \M_n(F(Z(R/P)))$ be all the possible options to view $\bar{x}-\lambda1$ as an element of $\M_n(Z(R/P))$.
For each $B_l$, $g(B_l)=0$, so $B_l$ is algebraic over $k$, hence (from arguments just seen for $A_j$), the characteristic polynomial of $B_l$, $c'_l(t)$, is over $k$ ($c'_l(t) \in k[t]$).

Notice that for any $A_j$, $A_j-\lambda I$ is a possible option for $\bar{x}-\lambda1$ as an element of $\M_n(Z(R/P))$, so $\{A_j-\lambda I\}_{j \in J} \subseteq \{B_l\}_{l \in L}$. Actually, one can see that $\{A_j-\lambda I\}_{j \in J} = \{B_l\}_{l \in L}$.

Choose any $j \in J$, then for $A_j-\lambda I$: $A_j-\lambda I=B_l$ for some $l \in L$.
Denote $c'_l(t)=(-1)^nt^n+b_{n-1}t^{n-1}+\ldots+b_1t+b_0$ with $b_i \in k$. Of course $b_0=det(B_l)$.
So, we have $b_0=det(B_l)=det(A_j-\lambda I)=c_j(\lambda) \in k-0$.
{}From Cayley-Hamilton Theorem, $0=c'_l(B_l)=(-1)^nB_l^n+b_{n-1}B_l^{n-1}+\ldots+b_1B_l+b_0$,
hence $(-1)^nB_l^n+b_{n-1}B_l^{n-1}+\ldots+b_1B_l=-b_0$, 
implying $((-1)^nB_l^{n-1}+b_{n-1}B_l^{n-2}+\ldots+b_1I)B_l=-b_0$,
so $(-b_0)^{-1}((-1)^nB_l^{n-1}+b_{n-1}B_l^{n-2}+\ldots+b_1I)B_l=I$, 

so $(-b_0)^{-1}((-1)^n(A_j-\lambda I)^{n-1}+b_{n-1}(A_j-\lambda I)^{n-2}+\ldots+b_1I)(A_j-\lambda I)=I$,

This shows that, in $R/P$, $(-b_0)^{-1}((-1)^n(\bar{x}-\lambda1)^{n-1}+b_{n-1}(\bar{x}-\lambda 1)^{n-2}+\ldots+b_11)(\bar{x}-\lambda1)=1$.

Therefore, in $R/P$, $\bar{x}-\lambda1$ is invertible (since we have found an element in $R/P$, $(-b_0)^{-1}((-1)^n(\bar{x}-\lambda1)^{n-1}+b_{n-1}(\bar{x}-\lambda 1)^{n-2}+\ldots+b_11)$ which is a left and a right inverse of $\bar{x}-\lambda1$).

We have shown that, for any $\lambda \in k-E_P$, $\bar{x}-\lambda1=x-\lambda +P$ is invertible in $R/P$ (we have already seen that $E_P \subset k$, $|E_P|\leq d$).
\end{proof}

Notice that in the above proof of Theorem \ref{images}, we viewed $x+P \in R/P$ as a matrix over $Z(R/P)$, although an embedding $R/P \hookrightarrow M_n(Z(R/P))$ does not necessarily exists.
In the next theorem we impose a similar condition to that of Theorem \ref{images}, and this time $R/P$ will be embedded in a matrix ring.

But first observe that for an affine prime Goldie non-simple non-PI $k$-algebra with GK dimension $<3$, where $k$ is an uncountable field, the following is true:

Take a prime ideal $P$ of $R$ with $\GKdim(R/P)=1$.
Since $R/P$ is a prime Goldie ring ($R/P$ is Goldie since it is prime PI), it has a classical (Goldie) ring of quotients, $Q(R/P)$, which is isomorphic to $\M_n(D)$, where $D$ is a division $k$-algebra (see \cite[Theorem 1.35]{gol}).
Actually, since $R/P$ is prime PI, Posner's theorem (\cite[Theorem 3.2]{gol}) says that $Q(R/P)$ is also PI, so by Kaplansky's theorem $Q(R/P)$ is finite dimensional over its center, $Z(Q(R/P))\cong Z(\M_n(D)) \cong Z(D)$.
So, $\dim_{Z(D)}(\M_n(D))<\infty$, hence $\dim_{Z(D)}(D)<\infty$.
Denote $\dim_{Z(D)}(D)=m$, then $\dim_{Z(D)}(\M_n(D))=mn^2$ and we have $\M_n(D)\hookrightarrow \End_{Z(D)}(\M_n(D))\cong \M_{mn^2}(Z(D))$.
Therefore, $R/P \hookrightarrow Q(R/P)\cong \M_n(D)\hookrightarrow \M_{mn^2}(Z(D))$.

$Z(D)$ is, of course, a field extension of $k$. We will add a condition on $Z(D)$, which, together with the existence of a special non-algebraic element, will guarantee the primitivity of $R$.
(Notice that Bell has mentioned in \cite[Remark 2.5]{bell} that $Z(D)$ is a finitely generated field extension of $k$ of transcendence degree 1. Hence, the condition we will add seems not too restrictive. Remember that when $k$ is algebraically closed, Theorem \ref{algebraic images} already tells us that $R$ is primitive).

\begin{thm}\label{images2}
Let $k$ be an uncountable field and let $R$ be an affine prime Goldie non-simple non-PI $k$-algebra with GK dimension $<3$. 
If $R$ has a non-algebraic element $x \in R$, such that $x+P \in R/P$ is algebraic over $k$, for every prime ideal $P$ with $\GKdim(R/P)=1$, and also for every such prime $P$, the center of the classical (Goldie) ring of quotients of $R/P$, $Z(Q(R/P))$, is isomorphic to $k(X)$, then $R$ is primitive.
\end{thm}

\begin{proof}
As was mentioned in the beginning of the proof of Theorem \ref{images}, it is enough to show that the following condition is satisfied: For every prime ideal $P$ with $\GKdim(R/P)=1$, there exists a finite (or countable) subset of $k$, call it $E_P$, such that for every $\lambda \in k-E_P$, $x-\lambda+P$ is invertible in $R/P$.

Let $P$ be a prime ideal with co-GK 1.
$\bar{x} \in R/P$ is algebraic over $k$; let $f(t) \in k[t]$ be such that $f(\bar{x})=\bar{0}$.
$\bar{x} \in R/P \hookrightarrow Q(R/P)\cong \M_n(D)\hookrightarrow \M_{mn^2}(Z(D))$, thus viewing $\bar{x}$ as an element of $\M_{mn^2}(Z(D))$, it has a characteristic polynomial $c(t) \in Z(D)[t]$ and a minimal polynomial $m(t) \in Z(D)[t]$.
Of course, in $Z(D)[t]$: $m(t)|f(t)$, so $f(t)=q(t)m(t)$ with $q(t) \in Z(D)[t]$.
But we have assumed that $Z(D)\cong k(X)$, so $f(t)=q(t)m(t)$ with $q(t),m(t) \in k(X)[t]$.
{}From Gauss's lemma, it follows that $q(t),m(t) \in k[X][t]$. But then, obviously, $q(t),m(t) \in k[t]$.
Now, in $k(X)[t]$, write $m(t)$ as a product of irreducible polynomials, each, apriori is over $k(X)$.
Then Gauss's lemma shows that each irreducible factor of $m(t)$ is over $k[X]$. But then each irreducible factor of $m(t)$ is over $k$. Therefore, $c(t)$ is also over $k$. Finally, take $E_P$ to be the set of roots of $c(t)$.
It is clear that for every $\lambda \in k-E_P$, $x-\lambda+P$ is invertible in $R/P$.
\end{proof}

\subsection{The extended center}
Smoktunowicz in her paper \cite{aga} proves the following:
"Let $k$ be a field, and let $R$ be an affine $k$-algebra, which is a domain with quadratic growth. Then either the center of $R$ is an affine $k$-algebra or R satisfies a polynomial identity (is PI) or else $R$ is algebraic over $k$".

We (almost) generalize Smoktunowicz's theorem in case the base field $k$ is uncountable.
\begin{thm}\label{fdfe}
Let $R$ be an affine prime Goldie non-simple $k$-algebra with GK dimension $<3$, where $k$ is an uncountable field.
Then $R$ is PI or the center of $R$ is a finite dimensional field extension of $k$.
\end{thm}
 
Observe that this is indeed a generalisation in case $k$ is uncountable:
a domain with quadratic growth is (prime) Goldie and quadratic growth implies GK dimension 2.
Also, in Theorem \ref{fdfe} (as well as in Smoktunowicz's theorem), non-simplicity of $R$ implies that $R$ is not algebraic over $k$. This is because an algebraic prime Goldie algebra is simple.

\begin{proof}
Assume that $R$ is non-PI. Hence $R$ is an almost PI algebra. Use Proposition \ref{centerfdfe} to conclude that $Z(R)$ is a finite dimensional field extension of $k$.
\end{proof}

We can consider the extended center of such algebras.

\begin{thm}\label{ecalg}
Let $R$ be an affine (non-simple) prime Goldie $k$-algebra with GK dimension $<3$, where $k$ is an uncountable field.
Then $R$ is PI or the extended center of $R$ is algebraic over $k$.
\end{thm}

\begin{proof}
Assume that $R$ is non-PI. Hence $R$ is an almost PI algebra. Use Proposition \ref{1.7.32} or Proposition \ref{extended} to conclude that $C(R)$ is algebraic over $k$.
\end{proof}

Theorem \ref{ecalg} is a special case of a result of Bell and Smoktunowicz brought in \cite{ec}: "Let $k$ be a field and let $A$ be a finitely generated prime $k$-algebra. We generalize a result of Smith and Zhang, showing that if $A$ is not PI and does not have a locally nilpotent ideal, then the extended centre of $A$ has transcendence degree at most $\GKdim(A) -2$ over $k$".
Indeed, our algebra $R$ is prime Goldie, hence $R$ does not have a locally nilpotent ideal. Therefore, their theorem says that the extended centre of $R$ has transcendence degree at most $\GKdim(R) -2 < 3-2=1$ over $k$, hence the extended centre of $R$ is algebraic over $k$.
Their theorem is more general because there $R$ is not necessarily Goldie, $\GKdim(R)$ can be $\geq 3$, and $k$ is not necessarily uncountable.

In the special case when $k$ is also algebraically closed, we get the following.

\begin{cor}\label{stably almost}
Let $R$ be an affine prime Goldie (non-simple) non PI $k$-algebra with GK dimension $<3$, where $k$ is an uncountable algebraically closed field.
Then $R$ is stably almost PI.
\end{cor}

\begin{proof}
Just use Corollary \ref{44} to conclude that $R$ is stably almost PI. 
\end{proof}

\section{Just infinite algebras}
\subsection{Some remarks about Farkas and Small's theorem}
Remember Farkas and Small's theorem which says that an affine semiprimitive just infinite algebra over an uncountable field is PI or primitive.

First notice that if one insists not to assume semiprimitivity, then one gets the following trichotomy (which is consistent with the Golod-Shafarevich example, which is algebraic).

\begin{thm}[Farkas and Small]\label{fs}
Let $R$ be an affine just infinite algebra over an uncountable field $k$.
Then $R$ is algebraic over $k$ or PI or primitive.
\end{thm}

\begin{proof}
Assume that $R$ in non-algebraic and non-PI.
Take any non-algebraic element $x \in R$. 
Denote all the nonzero prime ideals of $R$ by $\{Q_i\}_{i \in M}$.
Of course each nonzero prime ideal is maximal (since a prime finite dimensional algebra is simple).
So each nonzero prime ideal is of height 1 (remember that a just infinite algebra is prime, see \cite{prime}).
$R$ is an almost PI algebra (since we assume that $R$ is non-PI), hence from Proposition \ref{aa}, $|M|\leq \aleph_0$.

Define two countable subsets of $k$:
$A:=\{\alpha \in k | \ x-\alpha \text{ is left invertible in R}\}$: Indeed, from Amitsur's trick $A$ is countable ($x$ is non-algebraic).

$B:=\{\beta \in k | \ \text{there exists}\ i \in M \text{ such that} \ x-\beta+Q_i \in R/Q_i \\
\text{ is not invertible in } R/Q_i\}$: Indeed, $\forall i \in M$ $R/Q_i\hookrightarrow \End_k(R/Q_i)\cong \M_{\dim_k(R/Q_i)}(k)$.

Therefore, there exists only finitely many $\beta \in k$ (=the roots of the characteristic polynomial of $x+Q_i \in R/Q_i \hookrightarrow \End_k(R/Q_i)\cong \M_{\dim_k(R/Q_i)}(k)$) such that $x-\beta+Q_i=(x+Q_i)-(\beta+Q_i)$ is not invertible in $R/Q_i$, denote them $B(i)=\{\beta(i)_1,\ldots,\beta(i)_{m_i}\}$ (obviously, for $\beta \in k-B(i)$, $x-\beta +Q_i$ is invertible in $R/Q_i$, since its inverse, as an element of $\M_{\dim_k(R/Q_i)}(k)$, is a $k$-linear combination of powers of $x-\beta +Q_i$, which is in $R/Q_i$). Then $B=\cup_{i \in M}B(i)$, which is countable.

Next, define $E:=k-(A\cup B)$. $E$ is uncountable.

Take any $\lambda \in E$.

$x-\lambda$ is not left invertible in $R$ (since $\lambda \notin A$).
Hence $R(x-\lambda)\subsetneq R$, so there exists a maximal left ideal $L$ such that $R(x-\lambda)\subseteq L$.
$\ann_R(R/L)$ is a primitive ideal of $R$. 

There are two options for $\ann_R(R/L)$:
\bit
\item [(I)] $\ann_R(R/L)=0$. Then $R$ is primitive.
\item [(II)] $\ann_R(R/L) \neq 0$. Hence $\ann_R(R/L)=Q_i$ for some $i \in M$.
$Q_i=\ann_R(R/L)=\core(L)\subseteq L$.
We have an impossible situation:
On the one hand, $x-\lambda+Q_i$ is invertible in $R/Q_i$ (since $\lambda \notin B$).
On the other hand, $x-\lambda+Q_i\in L/Q_i\subseteq R/Q_i$, where $L/Q_i$ is a maximal left ideal of $R/Q_i$.
So $\ann_R(R/L) \in T$ cannot be true.
\eit
Hence only (I) is possible, so $R$ is primitive.
\end{proof}

We show that instead of $R$ affine over an uncountable field $k$, one can assume "the usual dimension-cardinality hypothesis", namely $dim_k(R)<|k|$.

\begin{thm}[Farkas and Small]\label{fs2}
Let $R$ be a just infinite algebra over a field $k$, with $dim_k(R)<|k|$.
Then $R$ is algebraic over $k$ or PI or primitive.
\end{thm}

\begin{proof}
Notice that in the above proof of Theorem \ref{fs}, we have used affinity twice (and uncountability of $k$ once):
\bit
\item When using Proposition \ref{aa}, to conclude that $R$, as an affine almost PI algebra (since we assume that $R$  is non-PI) has only countably many height 1 primes.
\item When using Amitsur's trick, to conclude that the set $A$ is countable. Not only affinity was used, but also uncountability of $k$.
\eit
However, under the weaker assumption that $dim_k(R)<|k|$, we can still have the same conclusions:
\bit
\item Instead of Proposition \ref{aa} we can use Proposition \ref{apACC}. Indeed, $R$ is an almost PI algebra (again we assume that it is non-PI), and as a just infinite algebra it satisfies ACC(ideals), so in particular, ACC(semiprimes).
\item Of course, in order to use Amitsur's trick, it is enough to have $dim_k(R)<|k|$.
\eit
Therefore, the proof of Theorem \ref{fs} is applicable here also, with the two new arguments just mentioned used, instead of the old ones.
\end{proof}

Notice that in Theorem \ref{fs} as well as in Theorem \ref{fs2}, semiprimitivity was not assumed, since we wanted to have a trichotomy.
However, if one wishes to have a dichotomy (namely, PI or primitive), then also "the usual dimension-cardinality hypothesis" can replace "affine over an uncountable field".
This answers what Farina and Pendergrass-Rice have written in \cite[section 3]{prime}, that it is unknown whether, except semiprimitivity, the other hypotheses can be relaxed. But maybe they have meant that it is unknown what happens over finite or countable fields; in that case the next theorem does not give an answer.

\begin{thm}
Let $R$ be a semiprimitive just infinite algebra over a field $k$, with $dim_k(R)<|k|$.
Then $R$ is PI or primitive.
\end{thm}

\begin{proof}
Exactly the same proof as the original proof of Farkas and Small will work here.
Just pay attention to the following:
\bit
\item As a non-PI just infinite algebra, $R$ has only countably many nonzero primitive ideals, as was explained there. Affinity is not needed.
\item Since $dim_k(R)<|k|$, Amitsur's linear independence trick can be used to show that $R$ is an algebraic algebra over $k$. 
\item Fisher and Snider's theorem can still be applied; Since $R$ is a non-PI just infinite algebra, it is an almost PI algebra. Hence Proposition \ref{ccsubset} can be used (indeed, in Proposition \ref{ccsubset}, affinity is not needed).
The difference is that instead of taking as a countable cofinal subset the set $\{T_n(R)\}$, we take the set $\{S_{n,d}\}$. 
\eit
\end{proof}

\subsection{When an affine prime Goldie algebra with GK dimension less then 3 over an uncountable field is just infinite}
Let $R$ be an affine prime Goldie non-simple non-PI $k$-algebra with GK dimension $<3$, where $k$ is an uncountable field.

It may happen that $R$ is a (non-PI) just infinite algebra. In that case, $R$ must be primitive:
Indeed, from Lemma \ref{2.semiprimitive} $R$ is semiprimitive, so Farkas and Small's theorem implies that $R$ is primitive.

If $R$ is bounded (for definition, see \ref{nilpotent-bounded}), then $R$ is a (non-PI) just infinite algebra. Indeed, we have seen in Theorem \ref{bounded} that a bounded affine almost PI algebra is (non-PI) just infinite.
However, boundedness is a very restrictive condition, so we look for a weaker condition that will guarantee that $R$ is just infinite.

\begin{prop}\label{3.countably}
Let $R$ be an affine prime Goldie non-simple non-PI $k$-algebra with GK dimension $<3$, where $k$ is an uncountable field.
If there are only countably many primitive ideals in $R$, then $R$ is a non-PI just infinite algebra.
In that case, $R$ is primitive.
\end{prop}

Notice that for the above $R$, the converse is also true, namely: If $R$ is a non-PI just infinite algebra, then there are only countably many primitive ideals in $R$. This was explained in the beginning of the proof of Farkas and Small's theorem.

\begin{proof}
In view of \cite[Corollary 1]{prime}, it is enough to show that every nonzero prime ideal has finite codimension.
{}From Lemma \ref{lemma}, $R$ can have three kinds of nonzero prime ideals.
If $R$ has a third kind prime, then Corollary \ref{2.cor} would imply that there are uncountably many first kind primes.
But each first kind prime is primitive, and we have assumed that $R$ has only countably many primitive ideals.
Therefore, $R$ has no third kind primes. Hence $R$ has no first kind primes (since a first kind prime must contain a third kind prime). So, $\Spec(R)=\{0,\text{second kind primes}\}$, and since each second kind prime has finite codimension, we are done.
\end{proof}

Smoktunowicz has written in \cite[page 266]{smo}:"...another question of Small: if $R$ is a Noetherian affine algebra with quadratic growth, does it follow that $R$ is either primitive or PI?
This is true in the graded case, as was shown by Artin and Stafford in 2000. According to Small, it is also true if every nonzero prime ideal in R is maximal." 

A remark about "According to Small, it is also true if every nonzero prime ideal in R is maximal":
\bit
\item If $k$ is uncountable: We shall see a proof of this in the following Proposition \ref{3.prop}.
\item If $|k|\leq \aleph_0$: It is interesting to see Small's proof. Apriori, it seems that such an algebra can even have a nonzero Jacobson radical, in contrast to the uncountable case where semiprimitivity is guaranteed (Lemma \ref{2.semiprimitive}).
\eit

In the same spirit, we bring the following Proposition \ref{3.prop}.
On the one hand, here we demand less: our algebra is Goldie with GK dimension $<3$, instead of Noetherian with quadratic growth and every nonzero prime ideal is primitive instead of maximal.
On the other hand, here we demand more: $k$ is uncountable, instead of any field and $R$ is prime.

\begin{prop}\label{3.prop}
Let $R$ be an affine prime Goldie non-simple non-PI $k$-algebra with GK dimension $<3$, where $k$ is an uncountable field.
If every nonzero prime ideal in $R$ is primitive, then $R$ is a non-PI just infinite algebra.
In that case, $R$ is primitive.
\end{prop}

\begin{proof}
First proof: $R$ is an affine almost PI algebra in which every nonzero prime ideal is primitive, hence Corollary \ref{1.9cor} says that $R$ is a non-PI just infinite algebra.

Second, more concrete, proof: From Lemma \ref{lemma}, $R$ can have three kinds of nonzero prime ideals.
We have assumed that every nonzero prime ideal in $R$ is primitive, so $R$ has no third kind primes. Hence $R$ has no first kind primes. So, $\Spec(R)=\{0,\text{second kind primes}\}$. We are done either by \cite[Corollary 1]{prime}, or by the following argument: there are only countably many second kind primes (remember that an affine almost PI algebra has only countably many height 1 primes, see Proposition \ref{aa}). Each second kind prime is primitive, so $R$ has only countably many primitive ideals. Now use Proposition \ref{3.countably}.
\end{proof}

Notice that for the above $R$, the converse is also true, namely: If $R$ is a non-PI just infinite algebra, then every nonzero prime ideal in $R$ is primitive (actually, maximal, by \cite[Corollary 1.6.30]{ro80}).

Therefore, for an affine prime Goldie non-simple non-PI $k$-algebra $R$ with GK dimension $<3$, where $k$ is an uncountable field: there are only countably many primitive ideals in $R$ $\iff$ $R$ is non-PI just infinite $\iff$ every nonzero prime ideal in $R$ is primitive.

\section{Acknowledgments}
I thank professor Louis Rowen and professor Uzi Vishne for the discussions we had about almost PI algebras.

\end{document}